\documentclass[twoside,11pt]{article}

\usepackage{amssymb, amsmath, latexsym, mathrsfs, verbatim, calc}

\usepackage[latin1]{inputenc}   

\def \N{{\rm I\!N}}

\newcommand{\dem}{\noindent \underline{Proof}: }
\newcommand{\cq}{\hfill$\Box$\\}
\newcommand{\ep}{\epsilon}

\def\i{{\bf 1}}

\def \R{{\mathop{{\rm I\negthinspace R}}}}

\def\qed{ \hfill \vrule width.25cm height.25cm depth0cm\medskip}
\def\Dt0{{\bf D}(t_0)}

\def\E{{\bf E}}
\def\M{{\bf M}}
\def\P{{\bf P}}

\def\p{{\bf p}}

\def\to{\rightarrow}

\def\O{\Omega}
\def\F{{\cal F}}
\def\D{\Delta}

\def\hp{{\hat p}}

\newcommand{\ba}{\[\begin{array}{rl}}
\newcommand{\ea}{\end{array}\]}
\newcommand{\be}{\begin{equation}}
\newcommand{\ee}{\end{equation}}
\newcommand{\bea}{\begin{eqnarray}}
\newcommand{\eea}{\end{eqnarray}}
\newcommand{\beaa}{\begin{eqnarray*}}
\newcommand{\eeaa}{\end{eqnarray*}}

\newtheorem{thm}{Theorem}[section]
\newtheorem{Lemma}[thm]{Lemma}
\newtheorem{cor}[thm]{Corollary}
\newtheorem{prop}[thm]{Proposition}
\newtheorem{rem}[thm]{Remark}
\newtheorem{eg}[thm]{Example}

\newcommand{\AR}{{\cal A}}

\newcommand{\FR}{{\cal F}}

\newcommand{\GR}{{\cal G}}
\newcommand{\HR}{{\cal H}}

\begin{document}

\title{On a continuous time game with incomplete information}

\author{Pierre Cardaliaguet\footnote{Laboratoire de Math\'ematiques, U.M.R CNRS 6205
Universit\'e de Brest,
6, avenue Victor-le-Gorgeu, B.P. 809, 29285 Brest cedex, France
e-mail : Pierre.Cardaliaguet@univ-brest.fr}
 and Catherine Rainer\footnote{Laboratoire de Math\'ematiques, U.M.R CNRS 6205
Universit\'e de Brest,
6, avenue Victor-le-Gorgeu, B.P. 809, 29285 Brest cedex, France
e-mail : Catherine.Rainer@univ-brest.fr}}

  \maketitle

\begin{abstract} 
For zero-sum two-player continuous-time games with integral payoff and incomplete information on one side, one shows that
the optimal strategy of the informed player can be computed through an auxiliary optimization problem over some  martingale measures. 
One also characterizes the optimal martingale measures and compute it explicitely in several examples. 
\end{abstract} 

\medskip
\noindent
{\bf \underline{Key words}:} Continuous-time game, incomplete information,  Hamilton-Jacobi equations.

\medskip
\noindent
{\bf \underline{MSC Subject classifications}: }  91A05, 49N70, 60G44, 60G42\\
$\;$\\

In this paper we investigate a two-player zero-sum continuous time game in which the players have an asymmetric information on the running payoff.
The description of the game involves
\begin{itemize}
\item[(i)] an initial time $t_0\geq 0$ and a terminal time $T>t_0$,
\item[(ii)] $I$ integral payoffs (where $I\geq 2$): $\ell_{i}:[0,T]\times U\times V \to\R$ for $i=1,\dots I$ where $U$ and $V$ are compact subsets of some finite dimensional spaces,
\item[(iii)] a probability $p=(p_i)_{i=1,\dots,I}$ belonging to the set $\Delta(I)$ of probabilities on $\{1, \dots, I\}$.
\end{itemize}
The game is played in two steps: at time $t_0$, the index $i$ is chosen at random among $\{1,\dots, I\}$ according
to the probability $p$~; the choice of $i$ is  communicated to Player 1 only.

Then the players choose their respective controls in order, for Player 1, to minimize the integral payoff
$\int_{t_0}^T\ell_{i}(s,u(s),v(s))ds$, and for Player~2 to maximize it. We assume that both players observe their opponent's control.
Note however that Player~2 does not  know which payoff he/she is actually maximizing.

Our game is a continuous time version of the famous repeated game with lack of information on one side studied by Aumann and Maschler
(see \cite{AM, So}). 
The existence of a value for our game has been investigated in \cite{cr} (in a more general framework):
if Isaacs' condition holds:
\begin{equation}\label{DefH}
H(t,p)=\inf_{u\in U}\sup_{v\in V} \sum_{i=1}^I p_i\ell_i(t,u,v)=\sup_{v\in V}\inf_{u\in U} \sum_{i=1}^I p_i\ell_i(t,u,v)\qquad \forall (t,p)\in [0,T]\times \Delta(I)\;,
\end{equation}
then the game has a value ${\bf V}={\bf V}(t_0,p)$ given by 
\be\label{valeur}
{\bf V}(t_0,p)=\inf_{(\alpha_i)\in ({\cal A}_r(t_0))^I} \sup_{\beta \in {\cal B}_r(t_0)}
\sum_{i=1}^I  p_i {\bf E}_{\alpha_i\beta}\left[\int_{t_0}^T \ell_{i}(s,\alpha_i(s), \beta(s))ds\right]
\ee
$$
\qquad =\sup_{\beta\in {\cal B}_r(t_0) }\inf_{(\alpha_i)\in ({\cal A}_r(t_0))^I}
\sum_{i=1}^I  p_i {\bf E}_{\alpha_i\beta}\left[\int_{t_0}^T \ell_{i}(s,\alpha_i(s), \beta(s))ds\right]\;,
$$
for any $(t_0,p)\in [0,T]\times \Delta(I)$,
where the $\alpha_i\in {\cal A}_r(t_0)$ (for $i=1,\dots,I$) are $I$ random strategies for Player 1, $\beta\in {\cal B}_r(t_0)$ is
a random strategy for Player 2 and ${\bf E}_{\alpha_i\beta}\left(\int_{t_0}^T \ell_{i}(s,\alpha_i(s), \beta(s))ds\right)$ is the payoff associated with the
pair of strategies $(\alpha_i,\beta)$: these notions are explained in the next section.
In \cite{cr} we also show that the value function ${\bf V}$ can be characterized in terms of {\it dual solutions} of some Hamilton-Jacobi
equations, which, following \cite{c5}, is equivalent to saying
that ${\bf V}$ is the unique viscosity solution of the following HJ equation:
\begin{equation}\label{min2}
\min\left\{ w_t+H(t,p)\,;\, \lambda_{\min}\left(\frac{\partial^2 w}{\partial p^2}\right)\right\}\;=\;0\qquad {\rm in }\; [0,T]\times \Delta(I)\;.
\end{equation}
In the above equation, $\lambda_{\min}\left(A\right)$ denotes the minimal eigenvalue of $A$, for any symmetric matrix $A$.
Note in particular that this equation says that ${\bf V}$ is convex with respect to $p$.\\

This paper is mainly devoted to the construction and the analysis of the optimal strategy for the informed player (Player 1). In particular we want to
understand how he/she has to quantify the amount information he/she reveals at each time. 
Our key step towards this aim is the following equality:
\be\label{mainrel}
{\bf V}(t_0, p_0)=\min_{{\bf P}\in \M(t_0,p_0)} {\bf E}_{{\bf P}}\left[ \int_{t_0}^T H(s,{\bf p}(s))ds \right] \qquad \forall (t_0,p_0)\in [0,T]\times \Delta(I)\;,
\ee
where $\M(t_0,p_0)$ is the set of martingale measures $\P$ on the set $\Dt0$ of c\`{a}dl\`{a}g processes $\p$ living in $\Delta(I)$ and such that $\p(t_0^-)=p_0$
and $\p(T)$ belong to the extremal points of $\Delta(I)$. Equality (\ref{mainrel}) is reminiscent of a result of \cite{DM}
in the discrete-time framework. It is directly related with the construction of the optimal strategy for the informed player. Indeed, let $\bar \P$ be
an optimal martingale measure in (\ref{mainrel}),  $\p(t)=(\p_1(t), \dots,\p_I(t))$ the coordonnate mapping 
 on $\Dt0$ and $\{e_1, \dots, e_I\}$ the canonical basis of $\R^I$. Then the informed player, knowing that the index $i$ has been chosen by nature, has just to play 
the random control $t\to {\rm argmin}_u \sum_j \p_j(t)\ell_j(t,u,v)$ with probability $\bar \P_i$, where $\bar \P_i$ is the restriction  of the law $\bar \P$ to the event
$\{\p(T)=e_i\}$. 

We give two proofs of (\ref{mainrel}). The first one relies
on a time discretization of  the value function taken from \cite{Sou} and the explicit construction of an approximated martingale measure in the discrete game.
The second one is based on a dynamic programming for the right-hand side of (\ref{mainrel}) and on the relation between this dynamic
programming and the notion of viscosity solution of (\ref{min2}). We shall see in particular that
the obstacle term $\lambda_{\min}\left(\frac{\partial^2 w}{\partial p^2}\right)\geq 0$ is directly related with the minimization over the martingale measures $\P$.

Since the optimal martingale in (\ref{mainrel}) plays a central role in our game, an analysis of this martingale is now in order. This
is the aim of section 3. We show that, under such a martingale measure 
$\bar \P$, the canonical process $\p(t)$ has to live in the set ${\cal H}(t)\subset \Delta(I)$ in which, heuristically, the Hamilton-Jacobi equation
$$
\frac{\partial {\bf V}}{\partial t}+H(t, p)=0
$$
holds. Moreover the process $\p(t)$ can only jump on the faces of the graph of ${\bf V}(t, \cdot)$. Namely, for any $t\in [t_0,T]$,
there is a measurable selection $\xi$ of $\partial {\bf V}(t, \p(t^-))$, such that
$$
{\bf V}(t,\p(t))-{\bf V}(t, \p(t^-)) -\langle \xi, \p(t)-\p(t^-)\rangle=0\qquad \bar \P-{\rm a.s.}
$$
Conversely, under suitable regularity assumptions on ${\bf V}$ and on ${\cal H}$, a martingale measure satisfying  these two conditions turns out to be  optimal in (\ref{mainrel}). \\

In section 4 we compute the optimal martingale measures for several examples of games. 
These examples show that, under the optimal measure, the process $\p$ can have very different 
behavior. When $I=2$ and under suitable regularity conditions, the optimal measure
is purely discontinuous. For instance we describe a game in which this optimal measure is unique and has to be the law of an Az\'ema martingale.
In higher dimension, uniqueness is lost in general. Moreover we give a class of examples in which, although there are optimal measures 
which are purely discontinuous,  they are also optimal measures under which $\p$ is a continuous process. 
For instance we show a game in which $\p$  lives on an expending
convex body moving with (reversed time) mean curvature motion. In this case case the law of $\p$ takes the form
$$
d\p_t=(I-a_t\otimes a_t)dB_t\;,
$$
where $(B_t)$ is a Brownian motion and the process $a_t$ lives in the unit sphere.

We complete the paper in section 5 by a list of open questions.

\section{Notations and general results}

In this part we introduce the main notations and assumptions needed in the paper. We also recall the main results of \cite{c1, cr}: existence
of a value for the game with lack of information on one side, as well as the characterization of the value.\\

{\bf Notations : } Throughout the paper, $x.y$ and $\langle x,y\rangle$ denote the scalar product in the space of vectors
$x,y\in \R^K$ (for some $K\geq 1$), and $|\cdot|$ the euclidean norm.  The closed ball of center $x$ and radius $r$ is denoted
by $B_r(x)$. The set $\Delta(I)$ is the set of probabilities measures on $\{1, \dots, I\}$, always identified with the simplex of $ R^I$:
$$
p=(p_1, \dots, p_I)\in \Delta(I)\quad  \Leftrightarrow \quad  \sum_{i=1}^Ip_i=1\; {\rm and }
\; p_i\geq 0\; {\rm for }\; i=1, \dots I\;.
$$
If $p\in \Delta(I)$, we denote by $T_{\Delta(I)}(p)$ the tangent cone to $\Delta(I)$ at $p$:
$$
T_{\Delta(I)}(p)=\overline{\bigcup_{\lambda>0} \left(\Delta(I)-p\right)/\lambda }\;.
$$
We denote by $\{e_i, \; i=1, \dots, I\}$ the canonical basis of $\R^I$. For any map $\phi:\Delta(I)\to \R$, ${\rm Vex}\phi$ is the convex  hull of $\phi$. 
We also denote by  ${\cal S}_I$ the set of symetric matrices of size $I\times I$. \\

Throughout the paper we assume the following conditions on the data:
\begin{equation}\label{HypfUV}
\left\{\begin{array}{rl}
i) & \mbox{\rm $U$ and $V$ are compact subsets of some finite dimensional spaces, }\\
ii) & \mbox{\rm For $i=1,\dots, I$, the payoff functions $\ell_i:U\times V\to \R$ are continuous.}
\end{array}\right.
\end{equation}
We also always assume that Isaac's condition holds and define the Hamiltonian of the game as:
\begin{equation}\label{Isaacs}
H(t,p):=\min_{u\in U}\max_{v\in V}\sum_{i=1}^I p_i\ell_i(t,u,v)=
\max_{v\in V}\min_{u\in U}\sum_{i=1}^I p_i\ell_i(t,u,v)
\end{equation}
for any $(t,p)\in  [0,T]\times  \Delta(I)$.\\

For any $t_0<T$, the set of open-loop controls for Player I is defined by
$$
{\cal U}(t_0)=\{u:[t_0,T]\to U \; \mbox{\rm Lebesgue measurable}\}\;.
$$
Open-loop controls for Player II are defined symmetrically and denoted ${\cal V}(t_0)$. 

A {\it pure} strategy for Player I at time $t_0$ is a map $\alpha: {\cal V}(t_0)\to {\cal U}(t_0)$  which is nonanticipative with delay, i.e.,
there is a partition $0\leq t_0<t_1<\dots<t_k=T$ such that,  for any $v_1,v_2\in {\cal V}(t_0)$, if $v_1= v_2$ a.e. on $[t_0,t_i]$ for some $i\in \{0, \dots, k-1\}$, then
$\alpha(v_1)= \alpha(v_2)$ a.e. on $[t_0,t_{i+1}]$. 

Let us fix ${\cal E}$ a set of probability spaces which is non trivial and stable by product.
A random control for Player I  at time $t_0$ is a pair $((\O_u, \F_u, \P_u), u)$ where the probability space $(\O_u, \F_u, \P_u)$ belongs to ${\cal E}$ and where
$u:\O_u \to {\cal U}(t_0)$ is Borel measurable from $(\O_u, \F_u)$ to ${\cal U}(t_0)$ endowed with the $L^1-$distance. We denote by
${\cal U}_r(t_0)$ the set of random controls for Player I  at time $t_0$ and abbreviate the notation $((\O_u, \F_u, \P_u), u)$ into simply $u$.

In the same way, a random strategy for Player I is a pair $((\O_\alpha, \F_\alpha, \P_\alpha), \alpha)$, where $(\O_\alpha,\F_\alpha,\P_\alpha)$ is a probability space in ${\cal E}$
and  $\alpha: \O_\alpha\times {\cal V}(t_0)\to {\cal U}(t_0)$ satisfies
\begin{itemize}
\item[(i)] $\alpha$ is measurable from $\O_\alpha \times {\cal V}(t_0)$ to ${\cal U}(t_0)$, with $\O_\alpha$ endowed with the $\sigma-$field $\F_\alpha$
and ${\cal U}(t_0)$ and ${\cal V}(t_0)$ with the Borel $\sigma-$field associated with the $L^1-$distance,

\item[(ii)] there is a partition $t_0<t_1<\dots<t_k=T$ such that,  for any $v_1,v_2\in {\cal V}(t_0)$, if $v_1\equiv v_2$ a.e. on $[t_0,t_i]$ for some $i\in \{0, \dots, k-1\}$, then
$\alpha(\omega, v_1)\equiv \alpha(\omega, v_2)$ a.e. on $[t_0,t_{i+1}]$ for any $\omega\in \O_\alpha$.
\end{itemize}
We denote by ${\cal A}(t_0)$ the set of pure strategies and by ${\cal A}_r(t_0)$ the set of random strategies for Player I. By abuse of notations,  an element of
${\cal A}_r(t_0)$ is simply noted  $\alpha$, instead of $((\O_\alpha, \F_\alpha, \P_\alpha), \alpha)$,
the underlying probability space being always denoted by $(\O_\alpha,\F_\alpha,\P_\alpha)$.
Note that ${\cal U}_r(t_0)\subset {\cal A}_r(t_0)$.

In order to take into account the fact that Player I knows the  index $i$ of the terminal payoff,
an admissible strategy for Player I is actually a $I-$uple $\hat{\alpha} =(\alpha_1,\dots,\alpha_I)\in ({\cal A}_r(t_0))^I$.\\

Pure and random controls and strategies for Player II are defined symmetrically; ${\cal V}_r(t_0)$
denotes the set of random controls for Player II, while ${\cal B}(t_0)$ (resp.  ${\cal B}_r(t_0)$)
denotes the set of pure strategies (resp. random strategies). Generic elements of ${\cal B}_r(t_0)$ are denoted by $\beta$, with
associated probability space $(\Omega_\beta,\F_\beta, \P_\beta)$.

Let us recall the
\begin{Lemma}[Lemma 2.2 of \cite{c1}]\label{PtFix} For any pair $(\alpha,\beta)\in{\cal A}_r(t_0)\times{\cal B}_r(t_0)$ and any $\omega:=(\omega_1, \omega_2)\in \O_\alpha\times \O_\beta$, there is a unique pair
$(u_\omega,v_\omega)\in{\cal U}(t_0)\times {\cal V}(t_0)$, such that
\begin{equation}\label{selle}
\alpha(\omega_1,v_\omega)=u_\omega \;{\rm and }\; \beta(\omega_2,u_\omega)=v_\omega\;.
\end{equation}
Furthermore the map $\omega\to (u_\omega,v_\omega)$ is measurable from $\O_\alpha\times \O_\beta$ endowed with $\F_\alpha\otimes \F_\beta$
into ${\cal U}(t_0)\times{\cal V }(t_0)$ endowed with the Borel $\sigma-$field associated with the $L^1-$distance.
\end{Lemma}

{\bf Notations : } Given any pair $(\alpha,\beta)\in{\cal A}_r(t_0)\times{\cal B}_r(t_0)$,
the expectation $\E_{\alpha\beta}$ is the integral
over $\O_\alpha\times \O_\beta$ against the probability measure $\P_\alpha\otimes\P_\beta$. In particular, if $\phi: [0,T]\times U\times V\to  \R$
is some bounded continuous map and $t\in (t_0,T]$,
we have
\begin{equation}\label{DefE}
{\bf E}_{\alpha\beta}\left[\int_t^T \phi(s,\alpha,\beta)ds\right]
:=
\int_{\O_\alpha\times \O_\beta} \left(\int_t^T \phi(s,u_\omega(s),v_\omega(s)ds\right) d\P_\alpha\otimes\P_\beta(\omega)\;,
\end{equation}
where $(u_\omega,v_\omega)$ is defined by (\ref{selle}). If one of the strategies is deterministic, we simply drop its subscript in 
the expectation. \\

As a particular case of Theorem 4.2 of \cite{cr} we have:

\begin{thm}[Existence of the value \cite{cr}]\label{value}  Assume that conditions (\ref{HypfUV}) and (\ref{Isaacs})
are satisfied. Then equality (\ref{valeur}) holds and we denote by ${\bf V}(t_0,p)$ the common value.
\end{thm}

In order to give the characterization of ${\bf V}$, we have to recall that the Fenchel
conjugate $w^*$ of a map $w:[0,T]\times \Delta(I)\to  \R$ is defined by
$$
w^*(t,\hp)=\max_{p\in \Delta(I)} p.\hp- w(t,p)\qquad \forall (t,\hp)\in [0,T]\times  \R^I\;.
$$
In particular ${\bf V}^{*}$ denotes the conjugate of ${\bf V}$.
If now $w$ is defined on the dual space $[0,T]\times \R^I$, we also denote by $w^*$ its conjugate with respect to $\hat{p}$ given by
$$
w^*(t,p)=\max_{\hp\in  R^I} p.\hp- w(t,\hp)\qquad \forall (t,p)\in [0,T]\times \Delta(I)\;.
$$

\begin{prop}[Characterization of the value, \cite{cr}]\label{CharVal} Under the assumptions of Theorem
\ref{value}, the value function ${\bf V}$ is the unique function defined on $[0,T]\times \Delta(I)$ such that
\begin{itemize}
\item[(i)] ${\bf V}$ is Lipschitz continuous in all its variables, convex with respect to $p$ and vanishes at $t=T$,

\item[(ii)]  for any $p\in \Delta(I)$, $t\to {\bf V}(t,p)$ is a viscosity subsolution of the Hamilton-Jacobi equation
\be\label{HJ}
w_t+H\left(t, p\right)=0\; {\rm in }\; (0,T)
\ee
where $H$ is defined by (\ref{Isaacs}),
\item[(iii)] For any smooth test function $\phi=\phi(\hp)$ and any $\hp\in \R^I$ such that $t\to {\bf V}^*(t,\hp)-\phi$ has a local maximum at
some $t$ and such that the derivative $p:=\frac{\partial {\bf V}^*}{\partial \hp}(t,\hp)$ exists, one has:
\be\label{HJDual}
\phi_t(t,\hp)-H\left(t, p\right)\geq 0\; .
\ee
\end{itemize}
We say that ${\bf V}$ is the unique {\it dual solution} of the HJ equation (\ref{HJ})
with terminal condition ${\bf V}(T,p)=0$.
\end{prop}

\begin{rem}{\rm In particular ${\bf V}$ does not depend on the class of probability space ${\cal E}$ chosen to define the random strategies. In view of the 
construction of the next chapter, we note that, given a family of probability measures, one can always built a set ${\cal E}$ containing this family and such that
${\cal E}$ is stable by product. 
}\end{rem}

In \cite{c5} we prove that the following characterization of ${\bf V}$ also holds:
\begin{prop}[Equivalent characterization \cite{c5}]\label{EquivChar} ${\bf V}$ is the unique Lipschitz continuous viscosity solution of the following obstacle problem
\begin{equation}\label{min3}
\min\left\{ w_t+H(t,p)\,;\, \lambda_{\min}\left(\frac{\partial^2 w}{\partial p^2}\right)\right\}\;=\;0 \;{\rm in }\; (0,T)\times \Delta(I)
\end{equation}
which satisfies ${\bf V}(T,p)=0$ in $\Delta(I)$.
\end{prop}

In the above proposition, we say that $w$ is a subsolution of the terminal time Hamilton-Jacobi equation
 (\ref{min3}) if, for any smooth test function $\phi: (0,T)\times \Delta(I)\to \R$ such that
$w-\phi$ has a local maximum at some point $(t,p)\in (0,T)\times Int(\Delta(I))$, one has
$$
\max\left\{ \phi_t(t,p)+H(t,p)\,;\, \lambda_{\min}\left(p,\frac{\partial^2 \phi}{\partial p^2}(t,p)\right)\right\}\;\geq \;0
$$
where, for any $(p,A)\in \Delta(I)\times {\cal S}_I$, 
$$
\lambda_{\min}(p, A)= \min_{z\in T_{\Delta(I)}(p)\backslash \{0\}}\langle Az,z\rangle/|z|^2\;.
$$We say that $w$ is a supersolution of (\ref{min3}) if, for any test function $\phi: (0,T)\times \Delta(I)\to \R$ such that
$w-\phi$ has a local minimum  at some point $(t,p)\in (0,T)\times\Delta(I)$, one has
$$
\max\left\{ \phi_t(t,p)+H(t,p)\,;\, \lambda_{\min}\left(p,\frac{\partial^2 \phi}{\partial p^2}(t,p)\right)\right\}\;\leq \;0\;.
$$
Finally a solution of (\ref{min3}) is a sub- and a super-solution of (\ref{min3}).

\section{Representation  of the solution}

Let us denote by $\Dt0$ the set of
c\`{a}dl\`{a}g functions from $\R\to \D(I)$ which are constant on $(-\infty, t_0)$
and on $[T,+\infty)$, by $t\mapsto {\bf p}(t)$ the coordinate mapping on $\Dt0$ and by ${\cal G}=({\cal G}_t)$ the filtration generated  by
$t\mapsto {\bf p}(t)$.

Given $p_0\in \Delta(I)$, we denote by $\M(t_0,p_0)$ the set of probability measures $\P$ on $\Dt0$ such that, under $\P$, $(\p(t),t\in[0,T])$ is a martingale and satisfies :
$$
{\rm for }\; t<t_0, \; \p(t)=p_0 \quad \mbox{\rm and, for }\;t\geq T,\;  {\bf p}(t)\in\{e_i,  \ i=1, \dots, I\}\qquad {\bf P}\ {\rm -a.s.}\;.
$$
Finally for any measure ${\bf P}$ on $\Dt0$, we denote by ${\bf E}_{{\bf P}}[\dots]$ the expectation with respect to ${\bf P}$.

Our main result is the following equality:

\begin{thm} \label{main}
\be\label{MartingalePb}
{\bf V}(t_0, p_0)=\inf_{{\bf P}\in \M(t_0,p_0)} {\bf E}_{{\bf P}}\left[ \int_{t_0}^T H(s,{\bf p}(s))ds \right] \qquad \forall (t_0,p_0)\in [0,T]\times \Delta(I)\;.
\ee
\end{thm}

We shall give two proofs of the result. The first one is based on a discretization procedure for ${\bf V}$, while
the second one uses a more direct approach of dynamic programming. Before this we show how to use the above
theorem to get optimal strategies for the first player.

\subsection{Construction of an optimal strategy}

We explain here how to use  Theorem \ref{main} to built an optimal strategy for the informed Player.
The construction of an optimal strategy for the non-informed Player, which uses completely different arguments (the so-called
approchability procedure) is described in \cite{Sou2}.

\begin{Lemma}\label{ExisOpti} For any $(t_0,p_0)$ there is at least one optimal martingale measure for problem
(\ref{MartingalePb}).
\end{Lemma}

\dem Let be a sequence of measures $(\P_n)_{n\in\N}\subset\M(t_0,p_0)$ satisfying
\begin{equation}
\label{optibis}
 {\bf V}(t_0,p_0)=\lim_{n\rightarrow+\infty}\E_{\P_n}\left[ \int_{t_0}^T H(s,{\bf p}(s))ds \right].
\end{equation}
Since, under all $\P\in\M(t_0,p_0)$, the coordinate process $\p$ is a martingale with support in the same compact space $\Delta(I)$, $(\P_n)$ converges weakly
(up to some subsequence)  to some measure $\bar \P$ that still belongs to $\M(t_0,p_0)$ (see Meyer-Zheng \cite{mezh}). Since $H$ is bounded and continuous,
passing to the limit in (\ref{optibis}) gives
$$
 {\bf V}(t_0,p_0)=\E_{\bar \P}\left[ \int_{t_0}^T H(s,{\bf p}(s))ds \right].
$$
Hence $\bar \P$ is optimal.
\qed

Let $(t_0,p_0) \in [0,T)\times \Delta(I)$ be fixed, $\bar \P$ be optimal in the problem (\ref{MartingalePb}).
Let us set $E_i=\{\p(T)=e_i\}$ and define the probability measure $\bar \P_i$ by:
$\forall A\in\GR$, $\bar \P_i(A):=\bar \P[A|E_i]=\frac{\bar \P(A\cap E_i)}{p_i}$, if $p_i>0$, and $\bar \P_i(A)=P(A)$ for an arbitrary probability measure $P\in\M(t_0,p_0)$ if $p_i=0$.\\
We also set
$$
\bar {\bf u}(t)=u^*(t,\p(t)) \qquad  \forall t\in \R
$$
and denote by $\bar {\bf u}_i$ the random control $\bar {\bf u}_i=((\Dt0, {\cal G}, \bar \P_i), \bar {\bf u})\in {\cal U}_r(t_0)$.

\begin{thm}\label{OptiStrat1} The strategy consisting in playing the random control $(\bar {\bf u}_i)_{i=1,\dots,I} \in ({\cal U}_r(t_0))^I$
is optimal for ${\bf V}(t_0,p_0)$. Namely
$$
{\bf V}(t_0,p_0)=\sup_{\beta\in {\cal B}_r(t_0)} \sum_{i=1}^I p_i{\bf E}_{\bar {\bf u}_i}\left[\int_{t_0}^T \ell_i(s,\bar {\bf u}_i(s), \beta(\bar {\bf u}_i)(s))ds\right]
$$
\end{thm}

\dem  Since
$$
\sup_{\beta\in {\cal B}_r(t_0)} \sum_{i=1}^I p_i{\bf E}_{\bar {\bf u}_i,\beta}\left[\int_{t_0}^T \ell_i(s,\bar {\bf u}_i(s), \beta(\bar {\bf u}_i)(s))ds\right]
=\sup_{\beta\in {\cal B}(t_0)} \sum_{i=1}^I p_i{\bf E}_{\bar {\bf u}_i}\left[\int_{t_0}^T \ell_i(s,\bar {\bf u}_i(s), \beta(\bar {\bf u}_i)(s))ds\right]\ ,
$$
it is enough to prove the equality
\begin{equation}\label{OptiStrateq}
{\bf V}(t_0,p_0)=\sup_{\beta\in {\cal B}(t_0)} \sum_{i=1}^I p_i{\bf E}_{\bar {\bf u}_i}\left[\int_{t_0}^T \ell_i(s,\bar {\bf u}_i(s), \beta(\bar {\bf u}_i)(s))ds\right]\;.
\end{equation}
Let us note that, since $\p(T)\in\{e_1, \dots, e_I\}$ $\bar \P$ a.s., we have
\be\label{toto1}
\sum_{i=1}^I {\bf 1}_{E_i}e_i=\p(T) \qquad \mbox{\rm $\bar \P$ a.s.}\;.
\ee
Let us fix a strategy $\beta\in {\cal B}(t_0)$ and set ${\bf v}=\beta(\bar{\bf u})$.
Then
$$
\sum_{i=1}^I p_i{\bf E}_{\bar {\bf u}_i}\left[\int_{t_0}^T \ell_i(s,\bar {\bf u}_i(s), \beta(\bar {\bf u}_i)(s))ds\right]
= \sum_{i=1}^I {\bf E}_{\bar \P}\left[\int_{t_0}^T \ell_i(s,\bar {\bf u}(s), {\bf v}(s))ds\ |\ E_i\right]{\bar \P}(E_i)
$$
$$
=   {\bf E}_{\bar \P}\left[\sum_{i=1}^I{\bf 1}_{E_i}\int_{t_0}^T \ell_i(s,\bar {\bf u}(s), {\bf v}(s)) ds\right]
$$
$$
= {\bf E}_{\bar \P}\left[\langle \p(T), \int_{t_0}^T \ell(s,\bar {\bf u}(s), {\bf v}(s)) ds\rangle \right]
$$
(from (\ref{toto1}) and where $\ell=(\ell_1, \dots, \ell_I)$ )
$$
= {\bf E}_{\bar \P}\left[ \int_{t_0}^T \langle \p(s),\ell(s,\bar {\bf u}(s), {\bf v}(s))\rangle ds \right]
$$
(by It\^{o}'s formula, since $\p(s)$ is a martingale under ${\bar \P}$ and the process $\ell(s,\bar {\bf u}(s), {\bf v}(s))$ is $({\cal G}_s)$ adapted)
$$
\leq {\bf E}_{\bar \P}\left[\int_{t_0}^T  \max_{v\in V} \langle \p(s),\ell(s,u^*(s,\p(s)), v)\rangle ds \right]
={\bf E}_{\bar \P}\left[\int_{t_0}^T  H(s, \p(s)) ds \right] ={\bf V}(t_0,p_0)\;.
$$
So we have proved that  (\ref{OptiStrateq}) holds.
\qed

\subsection{Proof of Theorem \ref{main} by discretization}\label{ProofDiscrete}

For simplicity of notations we shall prove Theorem \ref{main} for $t_0=0$, the proof in the general case being similar.

Let us start with an approximation of the value function. This approximation is a particular case of \cite{Sou}.
Let $n\in \N^*$, $\tau=1/n$ be the time-step and let us set $t_k=kT/n$ for $k=0,\dots, n$. We define by backward induction
$$
{\bf V}^\tau (T,p)=0\qquad \forall p\in \Delta(I)
$$
and, if ${\bf V}^\tau(t_{k+1},\cdot)$ is defined, then
$$
{\bf V}^\tau(t_{k},p)= {\rm Vex}\left( {\bf V}^\tau(t_{k+1},\cdot)+\tau H(t_k,\cdot)\right)(p)
$$
where ${\rm Vex}(\phi)(p)$ stands for the convex hull of the map $\phi=\phi(p)$ with respect to $p$.

\begin{Lemma}[\cite{Sou}] ${\bf V}^\tau$ uniformly converges to ${\bf V}$ as $\tau\to 0$ in the following sense:
$$
\lim_{\tau\to 0^+, \ t_k\to t, \ p'\to p} {\bf V}^\tau(t_k, p')={\bf V}(t,p)\qquad \forall (t,p)\in [0,T]\times \Delta(I)\;.
$$
\end{Lemma}

For any $k=0, \dots,n$ and $p\in \Delta(I)$ there are $\lambda^k=(\lambda^k_l)\in \Delta(I)$ and $\pi^k=(\pi^{k,l})\in \Delta(I)$ such that
\begin{equation}\label{defpi}
(i)\; \sum_{l=1}^{I}\lambda^k_l \pi^{k,l}=p\qquad {\rm and }\qquad
(ii)\; {\bf V}^\tau(t_{k},p)= \sum_{l=1}^{I} \lambda^k_l\left( {\bf V}^\tau(t_{k+1},\pi^{k,l})+\tau H(t_k,\pi^{k,l})\right)
\end{equation}
Without loss of generality we can choose the maps $p\mapsto \lambda^k(p)$ and $p\mapsto \pi^k(p)$ Borel measurable. For any $i\in \{1, \dots, I\}$
we now define a process $(\p^i_k, k\in\{ 0,\ldots, n+1\})$ on some arbitrary, big enough probability space $(\Omega,\FR,P)$ with values in $\Delta(I)$:\\
we start with $\p^i_0=p_0$ and, if $\p^i_k$ is defined for $k\leq n-1$, \\
\begin{itemize}
\item if the $i$-th coordinate of $\p^i_k$ satisfies $(\p^i_k)_i>0$, then the variable $\p^i_{k+1}$ takes its values in $\{\pi^{k,l}\left(\p^i_k\right), \; l\in 
\{ 1,\ldots, I\}\ \}$ with 
 \[ \forall l\in\{ 1,\ldots, I\},\; P[\p^i_{k+1}= \pi^{k,l}\left(\p^i_k\right)|\p^j_0,\ldots, \p^j_k,j=1,\ldots, I]=
\frac{\lambda^k_l(\p^i_k)\pi^{k,l}_i(\p^i_k)}{(\p^i_k)_i  },\]

\item if  $(\p^i_k)_i=0$, then we set $\p^i_{k+1}=\p^i_k$. 
\end{itemize}
For $k=n+1$, we  simply set
$
\p^i_{n+1}=e_i$ . \\
Finally we set $\p_k=\p^{{\bf i}}_k$ where ${\bf i}$ is the index chosen
at random by nature (i.e. ${\bf i}$ is a random variable that is independent from the processes $(\p^i_k), i\in\{1 ,\ldots, I\}$ and takes the values $1,\ldots,I$ with probability $p_1,\ldots,p_I$ respectively). The following Lemma is classical in the framework of repeated game theory with lack of information on one side
(see \cite{AM, So}).

\begin{Lemma}\label{marti}
 If we denote by $({\cal F}_k, k=0, \dots, n+1)$ the filtration generated by $(\p_k)$, then
$$
P\left[{\bf i}=i|{\cal F}_k\right]= (\p_k)_i\qquad \forall k\in \{0, \dots, n+1\},\ \forall i\in \{1,\dots, I\}\;.
$$
In particular, the process $(\p_k, k=0, \dots, n+1)$ is a martingale.
\end{Lemma}

\dem The result is obvious for $k=n+1$. Let us prove by induction on $k$ that, for all $i\in \{1,\dots, I\}$,
\begin{equation}
\label{marting}
P\left[{\bf i}=i|{\cal F}_{k}\right]= (\p_k)_i\; ,
\end{equation}
for $k=0,\ldots, n$. For $k=0$,  equality  (\ref{marting}) is obvious.
We now assume that (\ref{marting}) holds true up to some $k$ and check that  
 $P\left[{\bf i}=i \ |\ {\cal F}_{k+1}\right]=(\p_{k+1})_i$. 
 Since the variables $\p_k$ take their values in a finite set,   we can explicitely write
\begin{equation}
\label{i}
P\left[{\bf i}=i \ |\ {\cal F}_{k+1}\right]= \sum_{A\in\AR}
P\left[{\bf i}=i |A\right]
{\bf 1}_{A},
\end{equation}
where the  set $\AR$ is also finite and contains only sets of the form 
$A=\{ \p_1=\alpha_1,\ldots,\p_k=\alpha_k, \p_{k+1}=\pi^{k,l}(\p_{k})\}$ with $\alpha_1,\ldots,\alpha_k\in\Delta(I)$ and $l\in\{ 1,\ldots,I\}$ with $P[A>0]$. 
For such a $A\in\AR$, let us write 
\begin{equation}
\label{s1}
P[{\bf i}=i|A]=P[\{{\bf i}=i\}\cap A]/(\sum_{j=1}^IP[\{{\bf i}=j\}\cap A])\; .
\end{equation}
and,  for $j$ such that $(\alpha_k)_j>0$, using the independence between ${\bf i}$ and the processes $(\p_k^j)$,
\begin{equation}
\label{s2}\begin{array}{rl} 
P[\{{\bf i}=j\}\cap A]=
&P[{\bf i}=j,\p^j_1=\alpha_1,\ldots,\p^j_k=\alpha_k,\p_{k+1}^j=\pi^{k,l}(\p_{k}^j)]\\
=&P[{\bf i}=j]P[\p_{k+1}^j=\pi^{k,l}(\p_{k}^j)|\p^j_1=\alpha_1,\ldots,\p^j_k=\alpha_k]P[\p^j_1=\alpha_1,\ldots,\p^j_k=\alpha_k]\\
=&
P[{\bf i}=j,\p^j_1=\alpha_1,\ldots,\p^j_k=\alpha_k]\frac{\lambda^k_l(\alpha_k)\pi^{k,l}_j(\alpha_k)}{(\alpha_k)_j}.
\end{array}\end{equation}
Now, by the induction assumption,
\begin{equation}
\label{s3}\begin{array}{rl} 
P[{\bf i}=j,\p^j_1=\alpha_1,\ldots,\p^j_k=\alpha_k]=&
P[{\bf i}=j,\p_1=\alpha_1,\ldots,\p_k=\alpha_k]\\
=& (\alpha_k)_jP[\p_1=\alpha_1,\ldots,\p_k=\alpha_k]\; .
\end{array}\end{equation}
Thus, putting together (\ref{s1}),(\ref{s2}) and (\ref{s3}), we find out that
\[ P[{\bf i}=i|A]=\pi^{k,l}_i(\alpha_k).\]
But, on $A$, $\alpha_k=\p_k$, therefore, comming back to (\ref{i}), we get
\[ \begin{array}{rl}
P\left[{\bf i}=i \ |\ {\cal F}_{k+1}\right]= &\sum_{A\in\AR}
\pi^{k,l}_i(\p_k)
{\bf 1}_{A}\\
=&\sum_{l=1}^I\pi^{k,l}_i(\p_k){\bf 1}_{\{\p_{k+1}=\pi^{k,l}_i(\p_k)\}}\\
=&
(\p_{k+1})_i.
\end{array}\]
\qed

\begin{Lemma}\label{repVtau}
$$
{\bf V}^\tau(0,p_0)=E\left[ \tau \sum_{r=0}^{n-1} H(t_r,\p_{r+1})\right]\;.
$$
\end{Lemma}

\dem Let us show by (backward) induction on $k$ that
$$
E\left[ \tau \sum_{r=0}^{n-1} H(t_r,\p_{r+1})\right]=E\left[ {\bf V}^\tau(t_k,\p_k)+ \tau \sum_{r=0}^{k-1} H(t_r,\p_{r+1})\right]\;.
$$
Note that setting $k=0$ gives the Lemma.

For $k=n$, the result is obvious since ${\bf V}^\tau(T,p)=0$. Let us assume that the result holds true for $k+1$ and show that it still holds true for $k$.
By the induction assumption we have
$$
E\left[ \tau \sum_{r=0}^{n-1} H(t_r,\p_{r+1})\right]=E\left[ {\bf V}^\tau(t_{k+1},\p_{k+1})+ \tau \sum_{r=0}^{k} H(t_r,\p_{r+1})\right]\;.
$$
We note that, for all suitable function $f$,
\[\begin{array}{rl}
{\bf E}[f(\p_{k+1})|\sigma\{{\bf i},\p_k\}]=
&\sum_i\i_{\{ {\bf i}=i\}}\sum_l\frac{\lambda^k_l(\p^i_k)\pi^{k,l}_i(\p^i_k)}{(\p^i_k)_i}f(\pi^{k,l}(\p^i_k))\\
=&\sum_i\i_{\{ {\bf i}=i\}}\sum_l\frac{\lambda^k_l(\p_k)\pi^{k,l}_i(\p_k)}{(\p_k)_i}f(\pi^{k,l}(\p_k))
\end{array}\]
Thus, by Lemma  \ref{marti},
\[
{\bf E}[f(\p_{k+1})|\sigma\{\p_k\}]=\sum_i(\p_k)_i\sum_l\frac{\lambda^k_l(\p_k)\pi^{k,l}_i(\p_k)}{(\p_k)_i}f(\pi^{k,l}(\p_k))
=\sum_l\lambda^k_l(\p_k)f(\pi^{k,l}(\p_k)).
\]
In particular, by the definition of $\lambda^k$ and $\pi^k$ in (\ref{defpi}), we deduce that 
\[ {\bf E}[{\bf V}^\tau(t_{k+1},\p_{k+1})+\tau H(t_k,\p_{k+1})|\p_k]={\bf V}^\tau(t_k, \p_k).
\]
So
$$
E\left[ \tau \sum_{r=0}^{n-1} H(t_r,\p_{r+1})\right]=E\left[ {\bf V}^\tau(t_{k},\p_{k})+ \tau \sum_{r=0}^{k-1} H(t_r,\p_{r+1})\right]\;,
$$
which completes the proof.
\qed

We are now ready to show the inequality
\begin{equation}\label{ineq>}
{\bf V}(0, p_0)\geq \inf_{{\bf P}\in \M(0,p_0)} {\bf E}_{{\bf P}} \left( \int_{0}^T H(s,{\bf p}(s))ds \right)\;.
\end{equation}
Let $W(0,p_0)$ denote the right-hand side of this inequality.
Let us fix $\epsilon>0$  and let $\tau=1/n>0$ sufficiently small so that
$$
|{\bf V}(0, p_0)-{\bf V}^\tau(0,p_0)|\leq \epsilon
$$
and
$$
|H(t,p)-H(s,p)|\leq \epsilon \qquad \forall |s-t|\leq \tau, \; \forall p\in \Delta(I)\;.
$$
Let $(\p_k)$ be the martingale defined above. We built with this discrete time martingale a continuous one by setting:
$\tilde \p(t)=p_0$ if $t<0$, $\tilde \p(t)=\p_{k+1}$ if $t\in [t_k,t_{k+1})$ for $k=0, \dots, n-1$ and $\tilde \p(t)=\p_{n+1}$ if $t\geq T$.
The law of $(\tilde\p(t))$ defines a martingale measure ${\bf P}\in \M(0,p_0)$. Then
$$
\begin{array}{rl}
{\bf V}(0, p_0) \; \geq  & {\bf V}^\tau(0,p_0)-\epsilon\\
\geq & E\left[ \tau \sum_{r=0}^{n-1} H(t_r,\p_{r+1})\right]-\epsilon\\
\geq & E\left[\int_0^T H(s,\tilde \p(s))ds\right]-2\epsilon\\
\geq & {\bf E}_{{\bf P}}\left[\int_0^T H(s,\p(s))ds\right]-2\epsilon\\
\geq & W(0,p_0)-2\ep
\end{array}
$$
This proves (\ref{ineq>}).

Next we show
\begin{equation}\label{ineq<}
{\bf V}(0, p_0)\leq \inf_{{\bf P}\in \M(0,p_0)} {\bf E}_{{\bf P}} \left( \int_{0}^T H(s,{\bf p}(s))ds \right)\;.
\end{equation}
For this let us fix a martingale measure ${\bf P}\in \M(0,p_0)$. For $n\in \N^*$ large, we discretize the canonical process $\p$ in the usual way:
$\p^n(t)=\p(t_k)$ if $t\in [t_k,t_{k+1})$ for $k=0, \dots, n-1$ and $t_k=kT/n$. Then $\p^n(t)$ converges to $\p(t)$
$\P\otimes {\cal L}^1-$a.s. Therefore
$$
\lim_{n\to+\infty} {\bf E}_{{\bf P}}\left[\frac{1}{n} \sum_{r=0}^{n-1} H(t_r,\p^n(t_{r+1}))\right]= {\bf E}_{{\bf P}}\left[\int_0^T H(s,\p(s))ds\right]\;.
$$
To complete the proof of (\ref{ineq<}) it is enough to show  that
\begin{equation}\label{repVtaubis1}
{\bf E}_{{\bf P}}\left[\frac{1}{n} \sum_{r=0}^{n-1} H(t_r,\p^n(t_{r+1}))\right] \geq {\bf V}^\tau(0,p_0)
\end{equation}
where $\tau=1/n$. As for Lemma \ref{repVtau}, the proof of (\ref{repVtaubis1}) is achieved by showing by induction on $k\in \{0, \dots, n\}$ that
\begin{equation}\label{repVtaubis}
{\bf E}_{{\bf P}}\left[\frac{1}{n} \sum_{r=0}^{n-1} H(t_r,\p^n(t_{r+1}))\right] \geq {\bf E}_{{\bf P}}\left[{\bf V}^\tau(t_k,\p(t_k))+ \tau \sum_{r=0}^{k-1} H(t_r,\p^n(t_{r+1}))\right]
\end{equation}
For $k=n$ the result is obvious since ${\bf V}^\tau(T,\cdot)=0$.
Assume that the result holds for $k+1$.
We note that, since
$$
{\bf V}^\tau(t_{k+1},p)+\tau H(t_k,p) \geq {\bf V}^\tau(t_{k},p) \qquad \forall p\in \Delta(I)
$$
from the construction of ${\bf V}^\tau(t_k,\cdot)$ and since ${\bf V}^\tau(t_k,\cdot)$ is convex and $\p^n$ is a martingale under ${\bf P}$, we have
$$
\begin{array}{l}
{\bf E}_{{\bf P}}\left[ {\bf V}^\tau(t_{k+1},\p^n(t_{k+1}))+\tau H(t_k,\p(t_{k+1}))\ |\ {\cal G}_k\right]
\;  \geq \;   {\bf E}_{{\bf P}}\left[ {\bf V}^\tau(t_{k},\p^n(t_{k+1})) \ |\ {\cal G}_k\right] \\
\qquad  \geq \; {\bf V}^\tau\left(t_{k},{\bf E}_{{\bf P}}\left[\p^n(t_{k+1}) \ |\ {\cal G}_k\right] \right) \;  = \; {\bf V}^\tau(t_k, \p^n(t_k))
\end{array}
$$
So
$$
\begin{array}{rl}
{\bf E}_{{\bf P}}\left[\frac{1}{n} \sum_{r=0}^{n-1} H(t_r,\p^n(t_{r+1}))\right] \;  \geq & {\bf E}_{{\bf P}}\left[{\bf V}^\tau(t_{k+1},\p(t_{k+1}))+ \tau \sum_{r=0}^{k} H(t_r,\p^n(t_{r+1}))\right]\\
\geq &  {\bf E}_{{\bf P}}\left[{\bf V}^\tau(t_{k},\p(t_k))+ \tau \sum_{r=0}^{k-1} H(t_r,\p^n(t_{r+1}))\right]\;,
\end{array}
$$
which completes the proof of (\ref{repVtaubis}).

\qed

\subsection{Proof of Theorem \ref{main} by dynamic programming principle}

The idea is to show that the map
\begin{equation}
\label{opti}
W(t_0, p_0)=\inf_{{\bf P}\in \M(t_0,p_0)} {\bf E}_{\bf P}\left[ \int_{t_0}^T H(s,{\bf p}(s))ds \right] \qquad \forall (t_0,p_0)\in [0,T]\times \Delta(I)
\end{equation}
is a solution of equation (\ref{min3}) such that $W(T,p)=0$.
Since, in view of Proposition \ref{EquivChar} this equation has a unique solution and ${\bf V}$ is also solution, we get the desired result: $W={\bf V}$.

Arguing as in the proof of Lemma \ref{ExisOpti}, we first note that there is at least an optimal  martingale measure in the minimization problem (\ref{opti}).

\begin{Lemma}\label{lipsch}
$W$ is convex with respect to $p_0$ and  Lipschitz continuous in all variables.
\end{Lemma}

\dem Let $p_0,p_1 \in \Delta(I)$, $\P_0\in \M(t_0,p_0)$ and $\P_1\in \M(t_0,p_1)$ be optimal martingale measures
for the game starting from $(t_0,p_0)$ and $(t_0,p_1)$:
$$
{\bf E}_{\P_j}\left[ \int_{t_0}^T H(s,{\bf p}(s))ds \right] \;
= \; W(t_0,p_j)
$$
for $j=0,1$. For any $\lambda \in [0,1]$, let $p_\lambda=(1-\lambda)p_0+\lambda p_1$ and 
$\p_\lambda$ be the process on ${\bf D}(t_0)$ defined by 
\[ \p_\lambda(t)=\left\{\begin{array}{ll}
p_\lambda, &\mbox{ if } t< t_0\\
\p(t),&\mbox{ if } t\geq t_0.
\end{array}\right.\]
We define finally a probability measure $\P_\lambda\in\M(t_0,p_0)$ by :
For all measurable function $\Phi:{\bf D}(t_0)\rightarrow\R_+$,
\[ \E_{\P_\lambda}[\Phi(\p)]=(1-\lambda)\E_{\P_0}[\Phi(\p_\lambda)]+\lambda \E_{\P_1}[\Phi(\p_\lambda)].\]
 Then $\P_\lambda$
clearly belongs to $\M(t_0, p_\lambda)$ and
$$
W(t_0,p_\lambda) \leq {\bf E}_{\P_\lambda}\left[ \int_{t_0}^T H(t,{\bf p}(t))ds \right] = (1-\lambda)W(t_0,p_0)+\lambda W(t_0,p_1)\;.
$$
This proves that $W$ is convex with respect to $p_0$.

We now prove that $W$ is Lipschitz continuous. Since $W$ is convex with respect to $p_0$,  we need only to prove that $W$ is Lipschitz continuous
at the extremal points $(t_0,e_i),\; t_0\in[0,T],\;i\in\{ 1,\ldots, I\}$: namely we have to show that there is some $K\geq 0$ such that
$$
|W(t_0',p_0)-W(t_0,e_i)|\leq K( |p_0-e_i|+|t_0'-t_0|) \qquad \forall t_0',t_0\in [0,T], \; p_0\in  \Delta(I), \; i\in \{1, \dots, I\}\;.
$$
Let $t'_0, t_0\in [0,T]$, $p_0\in\Delta(I)$ and $i\in\{ 1,\ldots, I\}$.
One easily checks that $\M(t_0,e_i)$ consists in the single probability measure under which $(\p(s),t\leq s\leq T)$ is constant and equal to $e_i$. Consequently
$$
W(t_0,e_i)=\int_t^T H(s,e_i)ds\;.
$$
Let $\P$ be the optimal probability measure for problem (\ref{opti}) with starting point $(t_0',p_0)$.
We can write
\[\begin{array}{rl}
|W(t'_0,p_0)-W(t_0,e_i)|=&
|\E_\P\int_{t'_0}^TH(s,\p(s))ds-\int_{t_0}^T H(s,e_i)ds|\\
\leq & C|t_0-t'_0|+\E_\P\int_{t_0}^T|H(s,\p(s))-H(s,e_i)|ds\;,
\ea
where $C$ is an upper bound of the map $(t,p)\mapsto |H(t,p)|$ on $[0,T]\times\Delta(I)$.\\
Now let $\kappa$ be a Lipschitz constant for $H$. Then
\ba
\E_\P\int_{t_0}^T|H(s,\p(s))-H(s,e_i)|ds\leq &
\kappa\E_\P\int_{t_0}^T|\p(s)-e_i|ds\\
\leq &\sum_{j=1}^I\E_\P\int_{t_0}^T|\p_j(s)-\delta_{ij}|ds\\
= & \sum_{j\neq i}^I\E_\P\int_{t_0}^T\p_j(s)ds+\E_\P\int_{t_0}^T(1-\p_i(s))ds\\
= & (T-t_0)\sum_{j=1}^I|p_0^j-\delta_{ij}|\;,
\ea
where, for the last line, we used the fact that, under $\P$ and for all $j\in\{ 1,\ldots, I\}$, $\p_j$ is a martingale.
Finally
\[ \sum_{j=1}^I|p_0^j-\delta_{ij}|\leq \sqrt{I}| p_0-e_i|,\]
which completes the proof.
\qed

\begin{Lemma}\label{PgrDyn} The  following
dynamic programming holds: for any $\GR$-stopping time $\theta$ taking its values in $[t_0, T]$,
\begin{equation}\label{PgrDyn2}
W(t_0, p_0)= \inf_{\P\in  \M(t_0,p_0)}
{\bf E}_\P \left[ \int_{t_0}^{\theta} H(s,{\bf p}(s))ds +W(\theta, {\bf p}(\theta))\right] \ .
\end{equation}
\end{Lemma}

\dem
Let us introduce the subset $\M^f(t_0,p_0)$ of $\M(t_0,p_0)$ consisting in the martingale measures $\P$ on $\Dt0$ starting
from $p_0$ at time $t_0$ and  for which there is
a finite set $S\subset \Delta(I)$ such that any $\p\in {\rm Spt}(\P)$ satisfies $\p(t)\in S$ for $t\in [t_0,T]$ $\P$-a.s.
It is known that $\M^f(t_0,p_0)$ is dense in  $\M(t_0,p_0)$ for the weak* convergence of measures.
In particular it holds that
$$ W(t_0, p_0)=\inf_{\P\in  \M^f(t_0,p_0)} \E_\P\left[\int_{t_0}^TH(s,\p(s))ds\right] \;\; \forall (t_0,p_0)\in[0,T]\times\Delta(I),$$
and, since the map $\P\to {\bf E}_\P \left( \int_{t_0}^{\theta} H(s,{\bf p}(s))ds +W(\theta, {\bf p}(\theta))\right) $ is continuous for the weak* topology, the lemma is proved as soon we have shown that
\begin{equation}\label{infMf}
W(t_0, p_0)= \inf_{\P\in  \M^f(t_0,p_0)}
{\bf E}_\P \left[ \int_{t_0}^{\theta} H(s,{\bf p}(s))ds +W(\theta, {\bf p}(\theta))\right] .
\end{equation}
We shall prove (\ref{infMf}) for stopping times taking a finite number of values, then we generalize the result to all stopping times by passing to the limit.\\
 Let $\theta$ be a $\GR$-stopping time of the form
\begin{equation}
\label{esc}
 \theta=\sum_{l=1}^L\i_{A_l}\tau_l=\sum_{l=1}^L\i_{\{\theta=\tau_l\}}\tau_l,
\end{equation}
with $t_0\leq\tau_1<\ldots<\tau_L\leq T$ and, for all $l\in\{ 1,\ldots, L\}, A_l\in\GR_{\tau_l}$.
Let
 $\P\in\M^f(t_0,p_0)$ be $\ep-$optimal for $W(t_0,p_0)$
and $S=\{p^1, \ldots, p^K\}$ be such that $\P[\p(\theta)\in S]=1$.
We have
\begin{equation}
\label{dec2}
\E_{\P}\left[\int_{\theta}^T H(s,\p(s))ds\right]=\E_\P\left[
\sum_{j,l}\E_{\P}\left[\int_{\tau_l}^TH(s,\p(s))ds\Big|A_l\cap\{\p(\tau_l)=p^j\}\right]
\i_{A_l\cap\{\p(\tau_l)=p^j\}}\right]\; .
\end{equation}
But,  since $\P|_{A_l\cap\{\p(\tau_l)=p^j\}}\in\M(\tau_l,p^j)$, we have, for all $l$ and $j$,
\begin{equation}
\label{dec3}
\E_\P\left[\int_{\tau_l}^TH(s,\p(s))ds\Big|A_l\cap\{\p(\tau_l)=p^j\}\right] \geq W(\tau_l,p^j) \; .\end{equation}
Therefore
\begin{equation}
\label{unsens}
\begin{array}{rl}
W(t_0,p_0)+\epsilon\geq&\E_{\P}\left[\int_{t_0}^TH(s,\p(s))ds\right]\\
\geq &
\E_{\P}\left[\int_{t_0}^\theta H(s,\p(s))ds+
\sum_{l,j}W(\tau_l,p^j)\i_{ A_l\cap\{\p(\tau_l)=p^j\}}\right]\\
=& \E_{\P}\left[\int_{t_0}^\theta H(s,\p(s))ds+W(\theta,\p(\theta))\right]\; .
\end{array}\end{equation}

To prove the converse inequality, let $\P_0\in  \M^f(t_0,p_0)$ be $\ep-$optimal in the right-hand side of (\ref{infMf}) and let $S=\{p^1, \dots, p^K\}$ be such that $\P_0[\p(\theta)\in S]=1$.
For any  $\tau_l\in[t_0,T]$ and $p^j\in S$, let $\P_{l,j}$ be $\ep-$optimal for $W(\tau_l,p^j)$. 
We define a measure $\bar \P\in\M(t_0,p_0)$ in the following way :
For any function $f:\R\rightarrow \Delta(I)$ and $l\in\{ 1,\ldots, I\}$,
we define a process
\[ \p_{f,l}(t)=\left\{\begin{array}{ll}
f(t), &\mbox{ if } t< \tau_l\\
\p(t),&\mbox{ if } t\geq \tau_l.
\end{array}\right.\]
Then we set,
for all measurable function $\Phi:{\bf D}(t_0)\rightarrow\R_+$,
\[ \E_{\bar\P}[\Phi(\p)]=
\E_{\P_0}\left[\sum_{l,j}\i_{A_l\cap\{\p(\tau_l)=p^j\}}\E_{\P_{j,l}}[\Phi(\p_{f,l})]_{f=\p}\right].\]Then $\bar \P\in \M(t_0,p_0)$ and we have
\begin{equation}
\label{lotsens}
\begin{array}{rl}
W(t_0,p_0) \;  \leq & {\bf E}_{\bar \P} \left[\int_{t_0}^T H(s,{\bf p}(s))ds\right] \\
= & {\bf E}_{\P_0} \left[\int_{t_0}^{\theta} H(s,{\bf p}(s))ds+ \sum_{l,j} {\bf 1}_{A_l\cap\{\p(\tau_l)=p^j\}} {\bf E}_{\P_{l,j}} [\int_{\tau_l}^T H(s,{\bf p}(s))ds] \right] \\
\leq & {\bf E}_{\P_0} \left[\int_{t_0}^{\theta} H(s,{\bf p}(s))ds+ \sum_{l,j} {\bf 1}_{A_l\cap\{\p(\tau_l)=p^j\}} (W(\tau_l,p^j)+\ep) \right] \\
= & {\bf E}_{\P_0} \left[\int_{t_0}^{\theta} H(s,{\bf p}(s))ds+ W(\theta,\p(\theta))+\ep \right] \\
\leq & \inf_{\P\in  \M^f(t_0,p_0)}
{\bf E}_\P \left[ \int_{t_0}^{\theta} H(s,{\bf p}(s))ds +W(\theta, {\bf p}(\theta))\right]+2\ep\; .
\end{array}
\end{equation}
This allows us to conclude that (\ref{infMf}) holds for stopping times taking a finite number of values, and it remains now to show that (\ref{infMf}) holds for all stopping times.
But this last part of the proof is standard:  we just have to notice that,
if $\theta$ stands now for a general $\GR$-stopping time with $\theta\in [t_0,T]$, we can always find a sequence $(\theta_n)_{n\geq 0}$  of stopping times of the form (\ref{esc}) such that
$\theta_n \searrow\theta$ as $n\rightarrow\infty$ and that, for all $\P\in\M(t_0,p_0)$, we have
\[ \E_\P\left[\int_{t_0}^{\theta_n} H(s,\p(s))ds+W(\theta_n,\p(\theta_n))\right]\rightarrow_{n\rightarrow\infty}
\E_\P\left[\int_{t_0}^\theta H(s,\p(s))ds+W(\theta,\p(\theta))\right]\; .\]
 \qed

\begin{Lemma}\label{WSolVisc} $W$ is a solution of (\ref{min3}).
\end{Lemma}

\begin{rem} {\rm The proof of Lemma \ref{WSolVisc} is interesting because it shows that the martingale minimization problem gives rise to the penalization
term $\lambda_{\min}\left(p,\frac{\partial^2 \phi}{\partial p^2}\right)$ in (\ref{min3}).
}\end{rem}

\dem\ Let us first show that  $W$ is a supersolution of (\ref{min3}).   Let $\phi=\phi(t,p)$ be a smooth function such that $\phi\leq W$ with an equality at $(t_0,p_0)\in [0,T)\times \Delta(I)$.
We want to prove that
$$
\min\left\{ \phi_t(t_0,p_0)+H(t_0,p_0)\,;\,  \lambda_{\min}\left(p_0,\frac{\partial^2 \phi}{\partial p^2}(t_0,p_0)\right)\right\} \;\leq \;0\;.
$$
For this we assume that $\lambda_{\min}\left(p_0,\frac{\partial^2 \phi}{\partial p^2}(t_0,p_0)\right)>0$ and it remains to show that
\begin{equation}
\label{IneqSuper}
\phi_t(t_0,p_0)+H(t_0,p_0) \;\leq \;0\;.
\end{equation}

We claim that there are $r,\delta>0$ such that
\begin{equation}
\label{Afaire}
W(t,p)\geq \phi(t,p_0)+\langle \frac{\partial \phi}{\partial p}(t, p_0), p-p_0 \rangle+ \delta|p-p_0|^2 \qquad \forall t\in [t_0,t_0+r], \; \forall p\in \Delta(I)\;.
\end{equation}

{\it Proof of (\ref{Afaire}) : } From our assumption, there are $\eta>0$ and $\delta>0$ such that
$$
\langle \frac{\partial^2 \phi}{\partial p^2}(t, p)z,z\rangle  \geq 4\delta|z|^2  \qquad \forall z\in T_{\Delta(I)}(p_0), \ \forall (t,p)\in B_\eta(t_0,p_0)\;,
$$
where $T_{\Delta(I)}(p_0)$ is the tangent cone to $\Delta(I)$ at $p_0$:
$$
T_{\Delta(I)}(p_0)= \overline{ \bigcup_{h>0} (\Delta(I)-p_0)/h }
$$
Hence for $(t,p)\in B_\eta(t_0,p_0)$ we have
\begin{equation}\label{Afaire1}
W(t,p)\geq \phi(t,p)\geq \phi(t,p_0)+\langle \frac{\partial \phi}{\partial p}(t, p_0), p-p_0 \rangle+ 2\delta|p-p_0|^2\;.
\end{equation}
We also note that, for any $p\in \Delta(I)\backslash {\rm Int}(B_\eta(p_0))$, we have
\begin{equation}\label{Afaire2}
W(t_0,p)\geq \phi(t_0,p_0)+\langle \frac{\partial \phi}{\partial p}(t_0, p_0), p-p_0 \rangle+ 2\delta\eta^2
\end{equation}
because, if we set $p_1= p_0+\frac{p-p_0}{|p-p_0|}\eta$ and if $\hp_1\in \partial^-_pW(t_0,p_1)$, we have
$$
\begin{array}{rl}
W(t_0,p) \; \geq & W(t_0,p_1)+\langle \hp_1, p-p_1\rangle\\
\geq & \phi(t_0,p_0)+\langle \frac{\partial \phi}{\partial p}(t_0,p_0), p_1-p_0\rangle+ 2\delta\eta^2+\langle \hp_1, p-p_1\rangle\\
\geq & \phi(t_0,p_0)+\langle \frac{\partial \phi}{\partial p}(t_0,p_0), p-p_0\rangle+ 2\delta\eta^2+\langle \hp_1-
\frac{\partial \phi}{\partial p}(t_0,p_0), p-p_1\rangle
\end{array}
$$
where
$$
\langle \hp_1-
\frac{\partial \phi}{\partial p}(t_0,p_0), p-p_1\rangle
\geq 0
$$
because $w$ is convex, $\hp_1 \in  \partial^-_pW(t_0,p_1)$,  $\frac{\partial W}{\partial p}(t,p_0)\in\partial^-_p W(t_0,p_0)$ and $p-p_1=\gamma (p_1-p_0)$ for some $\gamma>0$.
Let us now argue by contradiction and assume that our claim (\ref{Afaire}) is false. Then there are $t_n\to t_0$
and $p_n\to p\in \Delta(I)$ such that
$$
W(t_n,p_n)< \phi(t_n,p_0)+\langle \frac{\partial \phi}{\partial p}(t_n,p_0), p_n-p_0 \rangle+ \delta|p_n-p_0|^2
$$
Note that $p_n\notin B_\eta(p_0)$ because of (\ref{Afaire1}).  Letting $n\to+\infty$, we get that $p\in \Delta(I)\backslash {\rm Int}(B_\eta(p_0))$ and
$$
W(t_0,p)\leq \phi(t_0,p)+\langle \frac{\partial \phi}{\partial p}(t_0,p_0), p-p_0 \rangle+ \delta\eta^2\;.
$$
This contradicts (\ref{Afaire2}). So (\ref{Afaire}) holds true for some $r>0$ sufficiently small. \\

Fix $\epsilon>0$ and $t\in (t_0,T)$.
Because of the dynamic programming (Lemma \ref{PgrDyn}), there exists $\P^t\in \M(t_0,p_0)$ such that
\begin{equation}
\label{utidyn}
{\bf E}_{\P^t} \left[\int_{t_0}^t H(s,{\bf p}(s))ds+ W(t, \p(t))  \right]
\leq W(t_0,p_0)+ \ep(t-t_0)
\end{equation}
Using the above inequality, (\ref{Afaire}) and the equality $\phi(t_0,p_0)=W(t_0,p_0)$ we get
\begin{equation}
\label{Hphidelta}
{\bf E}_{\P^t} \left[ \int_{t_0}^t H(s,{\bf p}(s))ds+ \phi(t, p_0)-\phi(t_0,p_0) + \delta|\p(t)-p_0|^2  \right]
\leq  \ep(t-t_0)
\end{equation}
because
$$
{\bf E}_{\P^t}  \left[ \langle \frac{\partial \phi}{\partial p}(t,  p_0), \p(t)-p_0 \rangle  \right] =0
$$
since $\P^t$ is a martingale measure. We note that (\ref{Hphidelta}) implies in particular that there is a constant $C>0$ such that
$$
{\bf E}_{\P^t} \left[\ |\p(t)-p_0|^2 \ \right]  \; \leq \; C(t-t_0)\qquad \forall t\in [t_0,t_0+r]
$$
because $H$ is bounded and $\phi$ is smooth. Since $({\bf p}(s))$ is a martingale under $\P^t$ this also
implies that
$$
{\bf E}_{\P^t} \left[\ |\p(s)-p_0|^2 \ \right] \; \leq \; C(t-t_0)\qquad \forall s, t\;{\rm with }\; t_0\leq s\leq t\leq t_0+r\;.
$$
Therefore, since $H$ is $\kappa-$Lipschitz continuous with respect to $p$, we have
$$
\left|{\bf E}_{\P^t} \left[ \int_{t_0}^t H(s,{\bf p}(s))ds\right]-\int_{t_0}^t H(s,p_0)ds\right|\leq
\kappa \int_{t_0}^t \left({\bf E}_{\P^t} \left[\ |{\bf p}(s)-p_0|^2\ \right]\right)^\frac12ds
\leq C\kappa(t-t_0)^\frac32\;.
$$
Plugging this inequality into (\ref{Hphidelta}) gives, for $t-t_0$ sufficiently small
$$
\int_{t_0}^t H(s,p_0)ds+ \phi(t, p_0)-\phi(t_0,p_0)\leq  2\ep(t-t_0)\;.
$$
Dividing this last inequality by $(t-t_0)$ and letting $t\to t_0$ gives (\ref{IneqSuper}).\\

Next we prove that $W$ is a subsolution of (\ref{min3}). This part relies on more classical arguments.
Let $\phi=\phi(t,p)$ be a smooth function such that $\phi\geq W$ with an equality at $(t_0,p_0)$
where $p_0\in Int(\Delta(I))$.  We want to prove that
$$
\min\left\{ \phi_t(t_0,p_0)+H(t_0,p_0)\,;\,  \lambda_{\min}\left(p_0, \frac{\partial^2 \phi}{\partial p^2}(t_0,p_0)\right)\right\} \;\geq \;0\;.
$$
We note that $\lambda_{\min}\left(p_0, \frac{\partial^2 \phi}{\partial p^2}(t_0,p_0)\right)\geq0$ because $W$ is convex with respect to $p$ and $p_0\in Int(\Delta(I))$.
So it remains to show that
$$
\phi_t(t_0,p_0)+H(t_0,p_0) \;\geq \;0\;.
$$
Fix $\epsilon>0$ and $t\in (t_0,T)$.
Because of the dynamic programming (Lemma \ref{PgrDyn}),
for any $\P\in \M(t_0,p_0)$ we have:
\begin{equation}
\label{utidyn2}
{\bf E}_\P \left[\int_{t_0}^t H(s,\p(s))ds+ W(t, \p(t))  \right]
\geq W(t_0, p_0)\;.
\end{equation}
Let us choose $\P=\delta_{p_0}$. Then we get from the definition of $\phi$:
$$
\int_{t_0}^t H(s,p_0)ds+ \phi(t, p_0)-\phi(t_0,p_0)
\geq  0\;.
$$
Dividing by $(t-t_0)$ and letting $t\to t_0$ gives the desired inequality since $\ep$ is arbitrary.
\qed

\noindent \underline{Proof of Theorem \ref{main} : } We have shown that $W$ is Lipschitz continuous (Lemma \ref{lipsch}),
that it  is a viscosity solution of equation (\ref{min3}) such that $W(T,p)=0$
(Lemma \ref{WSolVisc}). Since, from \cite{c5}, this equation has a unique Lipschitz continuous viscosity solution and since, from
Proposition \ref{EquivChar}, ${\bf V}$ is another Lipschitz continuous
viscosity solution of (\ref{min3}), we get the desired result: $W={\bf V}$.
\cq

As a consequence of Lemma \ref{PgrDyn} and of the above proof, we have:

\begin{cor}\label{PgrDynV} Let $(t_0,p_0)$ be an initial position and $\bar \P$ be an optimal martingale measure in (\ref{MartingalePb}).
Let $\theta\ge t_0$ be a stopping time. Then
$$
{\bf V}(t_0, p_0)= {\bf E}_{\bar \P} \left[ \int_{t_0}^{\theta} H(s,{\bf p}(s))ds +{\bf V}(\theta, {\bf p}(\theta))\right]\;.
$$
\end{cor}

\section{Analysis of the optimal martingale measure}

The section is devoted to the study of the optimal martingale measure in the optimization problem (\ref{MartingalePb}).
For doing so we first investigate the properties of the value function ${\bf V}$ as well as its conjugate ${\bf V}^*$. Then we define the set ${\cal H}\subset [0,T]\times \Delta(I)$
where---at least
heuristically--- ${\bf V}$ satisfies the Hamilton-Jacobi equation. We then show that, if $\bar\P$ is an optimal martingale measure, then the process
$\p$ remains on ${\cal H}$ and has jumps only on the flat parts of the graph of ${\bf V}(t,\cdot)$. These two conditions turn out to be sufficient to characterize
the optimal martingale measure under regularity assumptions on the value function.

\subsection{Some properties of ${\bf V}^*$}

We already know that ${\bf V}$ is a dual solution of the Hamilton-Jacobi equation (\ref{HJ}). In fact we have the  following sharper result:

\begin{prop}\label{EqV*} ${\bf V}^*$ is the solution of the Hamilton-Jacobi equation
\begin{equation}\label{HJ*}
\left\{\begin{array}{l}
\frac{\partial w}{\partial t} - H(t, \frac{\partial w}{\partial \hat p})=0  \qquad {\rm in }\; (0,T)\times \R^I  \\
w(T,\hat p)=\max\{\hat p_i\} \qquad {\rm in }\; \R^I
\end{array}\right.
\end{equation}
\end{prop}

\begin{rem}{\rm Compared to equation (\ref{HJDual}), where $\hat p$ appears as a parameter,  $\hat p$
is a genuine variable in the above equation. 
}\end{rem}

\dem We first show that ${\bf V}^*$ is a subsolution of (\ref{HJ*}).
Indeed, let $\phi$ be a ${\cal C}^1$ test function such that $\phi\geq {\bf V}^*$ with an equality at $(t,\hat p)$. Since ${\bf V}^*(t,\cdot)$ is convex and
$\phi$ is ${\cal C}^1$, this implies that ${\bf V}^*$ is differentiable with respect to $\hat p$ at $(t,\hat p)$ with $p=\frac{\partial {\bf V}^*}{\partial \hat p}(t,\hat p)=
\frac{\partial \phi}{\partial \hat p}(t,\hat p)$. Since $s\to {\bf V}^*(s,\hat p)-\phi(s,\hat p)$ has a maximum at $t$, the definition of dual solution gives
$$
\frac{\partial \phi}{\partial t}(t,\hat p)-H(t,p)=
\frac{\partial \phi}{\partial t}(t,\hat p)-H(t,\frac{\partial \phi}{\partial \hat p}(t,\hat p))\geq 0\;.
$$
So  ${\bf V}^*$ is a subsolution of (\ref{HJ*}).

Let $W$ be the solution of (\ref{HJ*}). Since ${\bf V}^*$ is a subsolution of this equation, we have $W\geq {\bf V}^*$ from the standard comparison
principle \cite{CIL}. In order to show the reverse inequality we are going to
check that $W^*$ is a dual supersolution of the Hamilton-Jacobi equation (\ref{HJ}).
 Since $H$ is positively homogeneous, independent of $\hat p$  and since $W(T,\cdot)$ is convex, $W(t,\cdot)$ is convex with respect to
$\hat p$ for any $t\in [0,T]$ (see \cite{GGIS}). From the usual representation formula for solutions of (\ref{HJ*}) (see \cite{ES}), we have
$$
W(t, \hat p)=\inf_{\beta \in {\cal B}(t)}\sup_{u\in {\cal U}(t)} \max\left\{ \hat p_i -\int_t^T \ell_i(s,\beta(u)(s),u(s))ds\right\}\;.
$$
We have
$$
W^*(t,p)= \sup_{\beta \in {\cal B}(t)} \sup_{\hat p\in \R^I} \min_{i\in I} \left\{ p.\hat p- \hat p_i +\inf_{u\in {\cal U}(t)} \int_t^T \ell_i(s,\beta(u)(s),u(s))ds\right\}
$$
If $p\in \Delta(I)$, for any $\beta\in {\cal B}(t)$ an optimum of the map
$$
\hat p\to  \min_{i\in I} \left\{ p.\hat p- \hat p_i +\inf_{u\in {\cal U}(t)} \int_t^T \ell_i(s,\beta(u)(s),u(s))ds\right\}
$$
is given by
$$
\hat p_i = \inf_{u\in {\cal U}(t)} \int_t^T \ell_i(s,\beta(u)(s),u(s))ds\; .
$$
Hence
$$
W^*(t,p)= \sup_{\beta \in {\cal B}(t)}\sum_{i=1}^I  p_i\inf_{u\in {\cal U}(t)} \int_t^T \ell_i(s,\beta(u)(s),u(s))ds\;,
$$
which is Lipschitz continuous in $p$. If $p\notin \Delta(I)$, a similar argument shows that $W^*(t,p)=+\infty$.
So the map $Z=W^*$ is Lipschitz continuous in all variables, convex in $p$, such that $Z^*=W$ is a subsolution of
the dual equation (\ref{HJ*}), which shows that $Z$ is a supersolution in the dual sense of Hamilton-Jacobi equation (\ref{HJ}). Since
${\bf V}$ is a dual solution, the comparison principle for dual solutions given in \cite{cr} implies that ${\bf V}\leq Z$, i.e., $W=Z^*\leq {\bf V}^*$. This shows that ${\bf V}^*=W$ is
the solution of (\ref{HJ*}).
\qed

\begin{prop}\label{differentiable} If $\partial {\bf V}^*(t,\hat p)=\{p\}$, then ${\bf V}^*$ is differentiable at $(t,\hat p)$ and
$$
\frac{\partial {\bf V}^*}{\partial t}(t,\hat p)=H(t,p)\;.
$$
\end{prop}

\dem Since $(s, \hat p')\to \partial {\bf V}(s,\hat p')$ is upper semicontinuous, for any $\epsilon>0$ one can find $\eta>0$
such that
$$
\partial {\bf V}(s,\hat p') \subset B_\epsilon(p) \qquad \forall (s,\hat p')\in {\cal O}:= (t-\eta,t+\eta)\times B_\eta(\hat p)\;.
$$
Hence  ${\bf V}^*$ satisfies
$$
\left| \frac{\partial {\bf V}^*}{\partial t}(s,\hat p') - H(t, p)\right|  \leq k\epsilon\;,
$$
for almost all $(s,\hat p')\in {\cal O}$, where $k$ is a Lipschitz constant of $H$. Thus
$$
\left| {\bf V}^*(s,\hat p')- {\bf V}^*(s',\hat p') - \int_s^{s'} H(\sigma, p)d\sigma\right| \leq  k\epsilon |s'-s| \qquad \forall (s,\hat p'),\  (s',\hat p) \in {\cal O}\;.
$$
Let now $(\tau_h,z_h)\to (\tau,z)$ in $\R\times \R^I$ as $h\to 0^+$. If $p_h \in \partial {\bf V}^*(t+h\tau_h, \hat p+hz_h)$, then $p_h\to p$
as $h\to 0$ and
$$
\begin{array}{rl}
{\bf V}^*(t+h\tau_h , \hat p+hz_h)-{\bf V}^*(t, \hat  p)
\;  \leq  &  {\bf V}^*(t, \hat p+hz_h)+\int_t^{t+h\tau_h} H(\sigma, p)d\sigma -{\bf V}^*(t, \hat  p)+\epsilon h|\tau_h|\\
\leq & h<p_h, z_h>+\int_t^{t+h\tau_h} H(\sigma, p)d\sigma +\epsilon h|\tau_h|
\end{array}
$$
Hence
$$
\limsup_{h\to 0^+} \frac{1}{h}\left({\bf V}^*(t+h\tau_h , \hat p+hz_h)-{\bf V}^*(t, \hat  p)
\right) \; \leq \; \langle p, z\rangle + H(t,p)\tau + \epsilon |\tau|
$$
In the same way one can prove that
$$
\liminf_{h\to 0^+} \frac{1}{h}\left({\bf V}^*(t+h\tau_h , \hat p+hz_h)-{\bf V}^*(t, \hat  p)
\right)\;  \geq \; <p, z>+ H(t,p)\tau - \epsilon |\tau|
$$
Since $\epsilon$ is arbitrary, we finally have the equality:
$$
\lim_{h\to 0^+} \frac{1}{h}\left({\bf V}^*(t+h\tau_h , \hat p+hz_h)-{\bf V}^*(t, \hat  p) \right) = <p, z>+ H(t,p)\tau \;,
$$
which shows that ${\bf V}^*$ is differentiable with $\frac{\partial {\bf V}^*}{\partial t}(t,\hat p)=H(t,p)$.
\qed


\subsection{The non revealing set ${\cal H}$ }

The aim of this section is the analysis of the set ${\cal H}$ defined by:
$$
{\cal H}= \left\{(t,p)\in [0,T)\times\Delta(I)\; |\;  \liminf_{h\to 0^+, \,  p'\to p}\frac{{\bf V}(t+h,p')-{\bf V}(t,p')}{h} = -H(t,p) \right\}\;.
$$
We also set
$$
{\cal H}(t)= \left\{p\in \Delta(I)\; |\; (t,p)\in {\cal H}\right\} \qquad \forall t\in [0,T)\;.
$$
In fact ${\cal H}$ is roughly speaking the set of points where  the Hamilton-Jacobi equation (\ref{HJ})
is satisfied. Indeed, if ${\bf V}$ is ${\cal C}^1$, then it is exactly so:
$$
{\cal H}= \left\{(t,p)\in [0,T)\times\Delta(I)\; |\;  \frac{\partial {\bf V}}{\partial t}+H(t,p)=0\right\}\;.
$$

\begin{Lemma}\label{IneqHJC} We have for any $(t,p)\in [0,T)\times \Delta(I)$,
\be\label{ine1}
{\bf V}(t,p)\leq {\bf V}(t+h, p)+\int_t^{t+h} H(\tau, p)d\tau \qquad \forall h\in [0,T-t]\;.
\ee
In particular,
\be\label{ine2}
\liminf_{h\to 0^+, \, t'\to t,  \, p'\to p}\frac{{\bf V}(t'+h,p')-{\bf V}(t',p')}{h} \geq -H(t,p)\;.
\ee
\end{Lemma}

\begin{rem} {\rm We have therefore:
$$
{\cal H}= \left\{(t,p)\in [0,T)\times\Delta(I)\; |\;  \liminf_{h\to 0^+, \,  p'\to p}\frac{{\bf V}(t+h,p')-{\bf V}(t,p')}{h} \leq -H(t,p) \right\}\;.
$$
}\end{rem}

\dem From the definition of dual solution, for any $p\in \Delta(I)$, the map $t\to {\bf V}(t,p)$ is a subsolution of HJ equation
(\ref{HJ}). Let $t<T$ and $p\in \Delta(I)$ be fixed. The solution $w$ of (\ref{HJ}) with terminal condition ${\bf V}(t+h,p)$  at time $t+h$
is given by the relation
$$
w(s)=\int_{s}^{t+h} H(\tau, p)d\tau +{\bf V}(t+h,p) \qquad \forall s\leq t+h\;.
$$
Since ${\bf V}(\cdot,p)$ is a solution of (\ref{HJ}) with the same terminal condition, we get
$$
{\bf V}(s,p)\leq w(s)=\int_{s}^{t+h} H(\tau, p)d\tau +{\bf V}(t+h,p) \qquad \forall s\leq t+h\;.
$$
Applying the above formula  to $s=t$ we get (\ref{ine1}). We then get (\ref{ine2}) thanks to the continuity of $H$.
\qed

\begin{Lemma} \label{HBorel} The set ${\cal H}$ is a Borel subset of $[0,T]\times \Delta(I)$ and ${\cal H}(t)$ is closed for any $t\in [0,T]$.
\end{Lemma}

\dem Indeed, ${\cal H}=\bigcap_{k}\bigcap_n {\cal H}_{nk}$ where
$$
{\cal H}_{nk}=\left\{(t,p)\in [0,T]\times  \Delta(I)\; |\;  \inf_{h\in (0,1/n], \,  |p'-p|\leq 1/n}\frac{{\bf V}(t+h,p')-{\bf V}(t,p')}{h} \leq -H(t,p) +\frac{1}{k} \right\}
$$
which are Borel subsets of $[0,T]\times  \Delta(I)$ since ${\bf V}$ and $H$ are continuous. Moreover ${\cal H}(t)$
is clearly closed for any $t\in [0,T]$ because $H$ is continuous.
\cq

In the proof of Lemma \ref{lipsch} we have already  noticed that
$$
{\bf V}(t,e_i)=\int_t^T H(s,e_i)ds\;,
$$
where $\{e_1,\ldots,e_I\}$ is the standard basis of $\R^I$.  This implies that $e_i\in\HR(t)$ for any $i=1, \dots, I$.
In particular, ${\cal H}(t)$ is nonempty for any $t\in [0,T]$. More precisely we have:

\begin{Lemma}\label{CScalH} Let $(t,p)\in [0,T)\times \Delta(I)$ and $\hp\in \partial {\bf V}(t,p)$. If ${\bf V}^*(t,\cdot)$ is differentiable at $\hp$, then
$(t,p)$ belongs to ${\cal H}$.
\end{Lemma}

\dem From Proposition \ref{differentiable}, ${\bf V}^*$ is differentiable at $(t,\hp)$ because ${\bf V}^*(t,\cdot)$ is differentiable at $\hp$.
 Let $h> 0$ be small and $p_h\in \partial {\bf V}^*(t+h, \hp)$. Then
$$
\frac{{\bf V}(t+h,p_h)-{\bf V}(t,p_h)}{h}\leq
- \frac{{\bf V}^*(t+h,\hp)-{\bf V}^*(t,\hp)}{h}
$$
because ${\bf V}(t+h,p_h)=p_h.\hp-{\bf V}^*(t+h,\hp)$ and ${\bf V}^*(t,\hp) \geq p_h.\hp-{\bf V}(t,p_h)$. Hence
$$
\limsup_{h\to0^+, \, p'\to p} \frac{{\bf V}(t+h,p_h)-{\bf V}(t,p_h)}{h}\leq
-\frac{\partial {\bf V}^*}{\partial t}(t,\hp) \leq -H(t,p)
$$
since ${\bf V}^*$ is differentiable at $(t,p)$ and from the definition of dual solution.
We complete the proof thanks to Lemma \ref{IneqHJC}
\qed

The next Lemma explains that, for any $(t,\hp)\in [0,T]\times \R^I$, the convex hull of $\partial {\bf V}^*(t,\hp)\cap {\cal H}(t)$ is exactly
equal to $\partial {\bf V}^*(t,\hp)$:

\begin{Lemma}\label{HCvHull} For any $(t,p)\in [0,T]\times \Delta(I)$ and any $\hp\in \partial {\bf V}(t,p)$, there are $(\lambda^j)\in \Delta(I)$,
$p^j\in {\cal H}(t)\cap \partial {\bf V}^*(t,\hp)$ for $j=1, \dots, I$ such that
$$
\sum_{j=1}^{I}\lambda^j p^j=p  \qquad {\rm and}\qquad  \sum_{j=1}^{I}\lambda^j {\bf V}(t,p^j)={\bf V}(t,p)\; .
$$
\end{Lemma}

\dem Let $\hp\in \partial {\bf V}(t,p)$ and, for $\ep>0$ small, $\hp_\ep=\hp+\ep\xi$. Since ${\bf V}^*$ is Lipschitz continuous, there are $\hp_n\to \hp_\ep$
at which ${\bf V}^*(t,\cdot)$ is differentiable. If we set $p_n=\frac{\partial {\bf V}^*}{\partial \hp}(t_n,\hp_n)$, then
Lemma \ref{CScalH} states that the points $p_n$ belong to ${\cal H}(t)$. Letting $n\to+\infty$, we can
find a subsequence of the $(p_n)$ which converges to some $p_\ep\in \partial {\bf V}^*(t,\hp_\ep)\cap  {\cal H}(t)$.

We now let $\ep\to 0$ to find some $p_\xi\in \partial {\bf V}^*(t,\hp)\cap {\cal H}(t)$. Moreover we have
$p_\xi.\xi\geq p.\xi$
because
$$
\langle p_\xi-p,\xi \rangle = \lim_{\ep\to0+}\langle p_\ep-p,\xi\rangle= \lim_{\ep\to0^+}\frac{1}{\ep}\langle p_\ep-p,\hp+\ep\xi-\hp\rangle \ \geq \ 0
$$
thanks to the monotony of the subdifferential. In particular,
$$
p.\xi\  \leq \ p_\xi.\xi \ \leq \ \sup_{p'\in \partial {\bf V}^*(t, \hp)\cap {\cal H}(t)} p'.\xi \qquad \forall \xi\in \R^N\;,
$$
which proves that $p$ belongs to the convex envelope of $\partial {\bf V}^*(t, \hp)\cap {\cal H}(t)$. The Lemma now follows easily
from the fact that ${\bf V}(t,\cdot)$ is affine on $\partial {\bf V}^*(t,\hp)$.
\qed

\subsection{Analysis of the optimal martingale measures}

We are now ready to study the optimal martingale measures in the optimization problem (\ref{MartingalePb}).
The main result of this section is the following:

\begin{thm}\label{bpH} Let ${\bar \P}$ be an optimal martingale measure in (\ref{MartingalePb}). Then
\be\label{pdansH}
(s,\p(s))\in {\cal H}\qquad \forall s\in [t_0,T],\; {\bar \P}\; {\rm a.s.}
\ee
and, for any $s\in (t_0,T]$, there is some measurable selection $\xi$ of $\partial {\bf V}(s, \p(s^-))$ such that
\be\label{Plats}
{\bf V}(s,\p(s))-{\bf V}(s,\p(s^-))-\langle \xi, \p(s)-\p(s^-)\rangle =0 \;\quad  {\bar \P}\; {\rm a.s.}\;.
\ee
\end{thm}

\begin{rem}{\rm The two conditions turn out to be necessary under suitable regularity conditions on the value function ${\bf V}$ and 
the martingale measure $\P$. See Theorem \ref{Verif} below.
}\end{rem}

\dem  For any $\ep, \delta>0$, let us set
$$
{\cal H}^c_{\ep,\delta}=\left\{(t,p)\in [0,T-\delta]\;|\; {\bf V}(t+h,p')-{\bf V}(t,p') \geq h(-H(t,p)+\ep) \; \forall (h,p')\in [0,\delta]\times B_\delta(p)\right\}\;.
$$
Then ${\cal H}^c_{\ep,\delta}$ is closed and
$$
\bigcup_{\ep,\delta>0} {\cal H}^c_{\ep,\delta}= ([0,T]\times\Delta(I))\backslash {\cal H}\;.
$$
Hence we have to prove that $(t,\p(t))\notin {\cal H}^c_{\ep,\delta}$ for any $t\in [t_0,T]$ $\bar \P-$a.s.
Let us note for later use that
\be\label{cloc}
{\bf V}(t+h,p)-{\bf V}(t,p) \geq - \int_t^{t+h} H(s,p)ds+\frac{\ep}{2} h \qquad \forall h\in [0,\delta], \; \forall (t,p)\in {\cal H}^c_{\ep,\delta}\;,
\ee
provided $\delta>0$ is small enough. Let us introduce the stopping time
$$
\theta=\inf\{ s\geq t \; |\; (s,\p(s))\in {\cal H}^c_{\ep,\delta} \}
$$
(with the convention that $\theta=T$ if $(s,\p(s))\notin {\cal H}^c_{\ep,\delta}$ for any $s\geq t$). Let $A=\{\theta<T\}$.
Let us assume that $\bar \P(A)>0$. From (\ref{cloc}) we have on $A$:
$$
{\bf V}(\theta, \p(\theta))\leq {\bf V}(\theta+h, \p(\theta))+\int_\theta^{\theta+h}H(s,\p(\theta))ds -\frac{\ep}{2}h\qquad \forall h\in [0,\delta]\;.
$$
Hence, for any $h\in [0,\delta]$,
$$
{\bf E}_{\bar \P}\left[{\bf V}(\theta, \p(\theta))\right] \leq
{\bf E}_{\bar \P}\left[ {\bf V}((\theta+h)\wedge T, \p(\theta))+\int_\theta^{(\theta+h)\wedge T}H(s,\p(\theta))ds\right] -\frac{\ep}{2}h\bar \P\left[A\right]\;.
$$
>From the dynamic programming principle Corollary \ref{PgrDynV} and the fact that $\bar \P$ is optimal we also have
$$
{\bf E}_{\bar \P}\left[{\bf V}(\theta, \p(\theta))\right] =
{\bf E}_{\bar \P}\left[ {\bf V}((\theta+h)\wedge T, \p(\theta))+\int_\theta^{(\theta+h)\wedge T}H(s,\p(s))ds\right]\;.
$$
So, for any $h\in (0,\delta]$, we have
$$
\frac{1}{h}{\bf E}_{\bar \P}\left[ \int_\theta^{(\theta+h)\wedge T}(H(s,\p(\theta))-H(s,\p(s))ds\right] \leq -\frac{\ep}{2}\bar \P\left[A\right]\;,
$$
which is impossible since $\p$ is right-continuous and $\bar \P[A]>0$. So we have proved that $\theta=T$ $\bar \P-$a.s., which means that
$(t,\p(t))\notin {\cal H}^c_{\ep,\delta}$ for any $t\in [t_0,T]$ $\bar \P-$a.s. \\

We now check that (\ref{Plats}) holds. Let $s>t_0$, $h>0$ and $\xi_h$ be a ${\cal G}_{s-h}$ measurable selection of $\partial {\bf V}(s-h, \p(s-h))$. Then we have
from the dynamic programming (Corollary \ref{PgrDynV})
$$
{\bf E}_{\bar \P}\left[ {\bf V}(s,\p(s))-{\bf V}(s-h,\p(s-h))-\int_{s-h}^s H(\tau, \p(\tau))d\tau\right]=0\;.
$$
Hence
\be \label{hpos}
\begin{array}{l}
{\bf E}_{\bar \P}\left[ {\bf V}(s,\p(s))-{\bf V}(s-h,\p(s-h))-\langle \xi_h, \p(s)-\p(s-h)\rangle\right]  \\
\qquad \qquad \leq\; {\bf E}_{\bar \P}\left[ \langle \xi_h, \p(s)-\p(s-h)\rangle\right] +h\|H\|_\infty=h\|H\|_\infty
\end{array}
\ee
since $\p$ is a martingale. Since $(\xi_h)$ is bounded in $L^\infty$, we can find a subsequence, again denoted $(\xi_h)$, which weakly
converges to some $\xi$ in $L^2$ as $h\to 0$. Note that $\xi\in \partial {\bf V}(s,\p(s^-))$ because $\xi_h\in \partial {\bf V}(s-h, \p(s-h))$
and $\p$ is has a left limit.
So we can let $h\to 0$ in (\ref{hpos}) to get
$$
 {\bf E}_{\bar \P}\left[ {\bf V}(s,\p(s))-{\bf V}(s,\p(s^-))-\langle \xi, \p(s)-\p(s^-)\rangle\right]\; \leq \; 0\;,
$$
where ${\bf V}(s,\p(s))-{\bf V}(s,\p(s^-))-\langle \xi, \p(s)-\p(s^-)\rangle\geq 0$ a.s. So (\ref{Plats}) holds.
\cq

\subsection{A verification Theorem}

If  the value function ${\bf V}$ is sufficiently smooth, then the conditions given in Theorem \ref{bpH} are ``almost sufficient"
in order to ensure a martingale measure to be optimal.

\begin{thm}\label{Verif} Let $(t_0,p_0)\in [0,T]\times \Delta(I)$.  Let us assume that ${\bf V}$ is of class ${\cal C}^{1,2}$  and that $\bar \P$ belongs to
$\M(t_0,p_0)$ and is such that
\begin{itemize}
\item[(i)] $\p(t)\in {\cal H}(t)$ for almost all $t\in [t_0,T]$ $\bar \P-$a.s.,

\item[(ii)] $\bar \P-$a.s.,
$$
{\bf V}(t,\p(t))-{\bf V}(t,\p(t^-))-\langle \frac{\partial {\bf V}}{\partial p}(t,\p(t^-)), \p(t)-\p(t^-)\rangle =0\qquad \forall t\in [t_0,T]\;,
$$

\item[(iii)] $\bar \P$ is a purely discontinuous martingale measure.
\end{itemize}
Then $\bar \P$ is optimal in problem (\ref{MartingalePb}).
\end{thm}

{\bf Remark : } The additional assumption that $\bar \P$ is purely discontinuous can be justified in some
particular cases. See Proposition \ref{CompDisc} below. \\

\dem Since ${\bf V}$ is of class ${\cal C}^{1,2}$, the set ${\cal H}$ is given by
$$
{\cal H}=\left\{(t,p)\in [0,T]\times \Delta(I)\; |\; \frac{\partial {\bf V}}{\partial t}(t,p)=-H(t,p)\right\}\;.
$$
We now use It\^{o}'s formula and the fact that $\bar \P$ is purely discontinuous to get
$$
\begin{array}{rl}
0={\bf E}_{\bar \P}\left[{\bf V}(T,\p(T))\right]= & {\bf V}(t_0,p_0)+{\bf E}_{\bar \P}\left[\int_{t_0}^T \frac{\partial {\bf V}}{\partial t}(s,\p(s))ds\right.\\
 & \left. +\sum_{s\geq t_0} \ {\bf V}(s,\p(s))-{\bf V}(s,\p(s^-))-\langle \frac{\partial {\bf V}}{\partial p}(s,\p(s^-)),\p(s)-\p(s^-)\rangle\right]
\end{array}
$$
$$
= {\bf V}(t_0,p_0)- {\bf E}_{\bar \P}\left[\int_{t_0}^T H(s,\p(s))ds\right]\;.
$$
The proof of Theorem \ref{Verif} is now complete thanks to Theorem \ref{main}.
\cq

\section{Examples}

\subsection{The autonomous case}

If the payoffs $\ell_i=\ell_i(u,v)$ are independent of time,  it is proved in \cite{Sou} that
\be\label{part1}
{\bf V}(t,p)=(T-t) {\rm Vex}H(p)\qquad \forall (t,p)\in [0,T]\times \Delta(I)\;.
\ee
Note that this equality is exactly what  Aumann-Maschler formula states for repeated games with incomplete information on one side (see \cite{AM}). In view of
(\ref{part1}) we have
$$
{\cal H}= [0,T]\times \left\{p\in\Delta(I)\; |\; {\rm Vex}H(p)=H(p)\right\}\;.
$$
Let us now fix $(t_0,p)\in [0,T]\times \Delta(I)$. Let $(\lambda_k)\in \Delta(I)$ and any $p^k\in \Delta(I)$ ($k=1, \dots, I\}$) such that
$$
\sum_{k=1}^I \lambda_k p^k=p\qquad {\rm and}\qquad \sum_{k=1}^I \lambda_k H(p_k)={\rm Vex}H(p)\;.
$$
We consider the probability measure $\bar\P\in\M(t_0,p)$ 
 under which, 
for all $k\in\{ 1,\ldots, I\}$, with probability 
$\lambda_k$, $\p$ is contant and equal to $p_k$ on $[t_0,T)$ .

\begin{prop} The measure measure $\bar \P$ is optimal for the minimization problem (\ref{MartingalePb}).
\end{prop}

\dem Indeed
$$
{\bf E}_{\bar \P}\left[\int_t^T H(\p(s))ds\right] = (T-t)\sum_{k=1}^I \lambda_k H(p^k)=(T-t){\rm Vex}H(p)={\bf V}(t,p)\;.
$$
\cq

\subsection{Examples when $I=2$}

In this section we assume that $I=2$.
We first show that, under suitable regularity properties of ${\bf V}$ and ${\cal H}$, there is a purely discontinuous martingale
which remains in ${\cal H}$ and jumps only on the flat parts of the graph of ${\bf V}$. Then we give an example where one can explicitly
compute the set ${\cal H}$ and the optimal martingale measures.

In this section we denote by $p\in [0,1]$ instead of $(p,1-p)$ (for $p\in [0,1]$) a generic element of $\Delta(I)$. The function ${\bf V}={\bf V}(t,p)$
will be defined on $[0,T]\times [0,1]$.

\begin{prop}\label{CompDisc} Let us assume that ${\bf V}$ is of class ${\cal C}^1$ and that
the set-valued map $t\to {\cal H}(t)$ enjoys the following regularity property:
there is some non decreasing map $K:[0,T]\to [0,+\infty)$ such that
\be \label{HLipschitz}
\forall s, t\in[0,T)\; {\rm with }\; s\leq t, \; \forall p\in {\cal H}(s), \; \exists p'\in {\cal H}(t) \; {\rm with }\; |p'-p|\leq K(t)-K(s)\;.
\ee
Then, for any initial position $(t_0,p_0)$ there is martingale measure $\bar \P\in \M(t_0,p_0)$ under which the process $\p$ satisfies conditions (i), (ii) and (iii) of Theorem \ref{Verif}.
\end{prop}
An example of value function satisfying condition (\ref{HLipschitz})  is given in Example \ref{ex1} below. 

\begin{rem}{\rm It is not known if there always exists an optimal martingale measure which is purely discontinuous without an additional assumption
like (\ref{HLipschitz}).  In fact in the case $I=2$ we have no example  of a martingale measure which satisfies (i) and (ii) of Theorem \ref{Verif} but not
(iii). For $I\geq 3$, we give an example below.
}\end{rem}

\dem Without loss of generality we assume that $t_0=0$. From Lemma \ref{HCvHull}, for any $s\leq t$,  and any $p\in {\cal H}(s)$, there are $p_1,p_2\in {\cal H}(t)$ such that
$p\in (p_1,p_2)$ and ${\bf V}(t,\cdot)$ is affine on $[p_1,p_2]$. If we choose $p_1$ as large as possible and $p_2$ as small as possible
(recall that we can do this since ${\cal H}(t)$ is closed), then we have from our assumption (\ref{HLipschitz}) that
$$
\min\{ |p-p_1|\ ; \  |p-p_2|\} \leq K(t)-K(s)\;.
$$
Let $\lambda\in [0,1]$ be such that $\lambda p_1+(1-\lambda)p_2=p$.  Note for later use that
\be\label{var}
\lambda |p_1-p|+(1-\lambda)|p_2-p|\leq 2(K(t)-K(s))
\ee
and that the maps $p_1=p_1(t,p)$, $p_2=p_2(t,p)$ and $\lambda=\lambda(t,p)$ Borel measurable.

Let us now introduce a large integer $n$ and a time step $\tau=T/n>0$. We set $t_k=\tau k$ for $i=0, \dots, n$. As in section \ref{ProofDiscrete} we define by induction the
process $(\p^n_k)_{k=-1, \dots, n}$ such that\\
(i) $\p^n_{-1}=p_0$, \\
(ii) for any $k\ge 0$, $\p^n_k\in {\cal H}(t_k)$, \\
(iii) knowing $\p^n_k$, $\p^n_{k+1}$ is equal to $p_1(t_k,\p^n_{1})$ with probability $\lambda(t_k,\p^n_{1})$ and $p_2(t_k,\p^n_{1})$ with probability $(1-\lambda(t_k,\p^n_{1}))$.\\
We first note that $\p^n$ is a martingale. From (\ref{var}) we have
$$
{\bf E}\left[ |\p^n_{k+1}-\p^n_k|\; |\; \p^n_k\right]\leq 2(K(t_k))-K(t_{k+1}))\;.
$$
Therefore the process $\p^n$ has bounded total variations:
\be\label{VTestimate}
{\bf E}\left[ \sum_{k=0}^{n-1}|\p^n_{k+1}-\p^n_k|\right]\leq 2(K(T)-K(0))\;.
\ee
We now interpolate the process $\p^n$ as in section \ref{ProofDiscrete}
 in order to get a martingale measure $\P^n\in \M(0,p_0)$.
Following \cite{mezh}, letting $n\to +\infty$, we can find a subsequence, again denoted $\p^n$, such that the law of the process $\p^n$ converges to some
$\bar \P\in M(0,p_0)$
and such that $\p^n(t)$ converges in law to $\p(t)$ for any $t$ belonging to some subset of full measure ${\cal T}$ of $[0,T]$.
Because of (\ref{VTestimate}), $\p$ has finite total variations under $\bar \P$ and therefore is purely discontinuous.
We now check that $\p$ satisfies conditions (i) and (ii).  Let ${\cal T}_1$ be the set of $t\in {\cal T}$ at which
the map $K$ is continuous. Then ${\cal T}_1$ is of full measure in $[0,T]$.
For any $t\in {\cal T}_1$ let $k_n$ be such that $k_n\tau\to t$ and $t\in [k_n\tau, (k_n+1)\tau)$.
   From assumption (\ref{HLipschitz}),
$$
d(p, {\cal H}(t)) \leq K(t)-K(t_{k_n}) \qquad \forall p\in {\rm Spt}(\p^n(t))\;,
$$
(where $d(p, {\cal H}(t))$ is the distance of $p$ to the set ${\cal H}(t)$)
because $\p^n(t)=\p^n_{k_n}$ and $\p^n_{k_n}\in {\cal H}(t_k)$ $\P-$a.s..
Letting $n\to +\infty$ implies that ${\rm Spt}(\p(t)) \subset {\cal H}(t)$ $\bar \P-$a.s. since $\p^n(t)$ converges in law to $\p(t)$ and $K$ is continuous at $t$.

For proving that $\p$ satisfies (ii), let us first note that
$$
{\bf V}(t_{k+1}, \p^n_{k+1})-
{\bf V}(t_{k+1}, \p^n_{k})- \langle \frac{\partial {\bf V}}{\partial p}(t_{k+1}, \p^n_{k}), \p^n_{k+1}-\p^n_{k}\rangle =0\qquad \forall k\in \{0, \dots, n-1\} \;.
$$
Hence
\be\label{EVk+1}
{\bf E}\left[{\bf V}(t_{k+1}, \p^n_{k+1})\right]={\bf E}\left[{\bf V}(t_{k+1}, \p^n_{k})\right]\qquad \forall k\in \{0, \dots, n-1\}
\ee
because $(\p^n_k)_k$ is a martingale.
Let $s,t\in {\cal T}$ be such that $s<t$ and $k_1, k_2$ be such that  $s\in [t_{k_1}, t_{k_1+1})$, $t\in [t_{k_2}, t_{k_2+1})$. Then
$$
{\bf V}(t,\p^n(t))-{\bf V}(s,\p^n(s)) \leq
{\bf V}(t_{k_2}, \p^n_{t_{k_2}})-{\bf V}(t_{k_1}, \p^n_{t_{k_1}})+ 2C(t-t_{k_2}+s-t_{k_1})
$$
where $C=\|\frac{\partial {\bf V}}{\partial t}\|_\infty$ and where
$$
\begin{array}{rl}
{\bf V}(t_{k_2}, \p^n_{t_{k_2}})-{\bf V}(t_{k_1}, \p^n_{t_{k_1}}) \;
= & \sum_{l=k_1}^{k_2-1} \left( \ {\bf V}(t_{l+1}, \p^n_{t_{l+1}})-{\bf V}(t_{l}, \p^n_{t_{l}}) \ \right)  \\
\leq & \sum_{l=k_1}^{k_2-1} \left( \ {\bf V}(t_{l+1}, \p^n_{t_{l+1}})-{\bf V}(t_{l+1}, \p^n_{t_{l}})\  \right) + C(t_{k_2}-t_{k_1})
\end{array}
$$
Combining (\ref{EVk+1}) with the above inequality gives
$$
{\bf E}\left[{\bf V}(t,\p^n(t))-{\bf V}(s,\p^n(s))\right] \leq C(t-s+\tau)\;.
$$
Letting $n\to +\infty$ leads to
$$
{\bf E}_{\bar \P}\left[\ {\bf V}(t,\p(t))-{\bf V}(s,\p(s)) \ \right] \leq C(t-s)
$$
>From the right-continuity of the process $\p(t)$, this inequality also holds for any $t$. Since $\p$ is a martingale we get
$$
{\bf E}_{\bar \P}\left[{\bf V}(t,\p(t))-{\bf V}(s,\p(s)   -\langle \frac{\partial {\bf V}}{\partial p}(t, \p(s)), \p(t)-\p(s)  \rangle             )\right] \leq C(t-s)
$$
for any $t\in (0,T]$,  $s\in {\cal T}$,  $s<t$. Letting now $s\to t^-$ with $s\in{\cal T}$ gives
$$
{\bf E}_{\bar \P}\left[{\bf V}(t,\p(t))-{\bf V}(t,\p(t^-)) -\langle \frac{\partial {\bf V}}{\partial p}(t, \p(t^-)), \p(t)-\p(t^-)  \rangle             )\right] \leq 0 \qquad \forall t\in (0,T]\;.
$$
Since ${\bf V}(t, \cdot)$ is convex this last inequality finally implies that
$$
{\bf V}(t,\p(t))-{\bf V}(t,\p(t^-)) -\langle \frac{\partial {\bf V}}{\partial p}(t, \p(s)), \p(t)-\p(s)  \rangle=0
$$
for any $t\in (0,T]$ $\bar P$ a.s.

\cq

Our aim is to identify ${\cal H}$ under the following assumption on $H$:
\begin{eg}\label{ex1} We assume that
there exist  $h_1,h_2:[0, T]\to [0,1]$ continuous, $h_1\leq h_2$, $h_1$ decreasing and $h_2$ increasing,   such that
\be\label{DefH1}
\mbox{\rm ${\rm Vex}H(t,p)=H(t,p)$ $\; \Leftrightarrow\; $   $p\in [0, h_1(t)]\cup [h_2(t),1]$ }
\ee
and
\be\label{DefH2}
\frac{\partial^2 H}{\partial p^2}(t,p) >0 \qquad \forall (t,p) \; {\rm with}\;  p\in [0, h_1(t))\cup (h_2(t),1]\;.
\ee
\end{eg}

\noindent {\bf For instance, } if we assume that $U=[-1,1]$, $V= [0,2\pi]$ and
$$
\ell_1(t,u,v)=u+ \alpha(t)\cos(v), \; \ell_2(t,u,v)=-u+ \alpha(t)\sin(v)\qquad \forall (u,v)\in U\times V
$$
where the smooth map $\alpha:[0,T]\to \R$ is decreasing and such that $\alpha(t)>2$ for any $t\in [0,T]$, then
$$
H(t,p) =-|2p-1|+ \alpha(t)\sqrt{p^2+(1-p)^2}
$$
satisfies (\ref{DefH1}) and (\ref{DefH2}) with $h_1(t)=1/2 - 1/(2\alpha^2(t)-4)^\frac12$, $h_2(t)=1/2 + 1/(2\alpha^2(t)-4)^\frac12$.

\begin{prop}\label{ex1H} Under the assumptions of Example \ref{ex1},
$$
{\bf V}(t,p)=\int_t^T {\rm Vex}H(s,p)ds \qquad \forall (t,p)\in [0,T]\times \Delta(I)
$$
and
\be\label{HpartCase}
{\cal H}=\left\{(t,p)\in [0,T]\times [0,1]\; |\; p\in[0,h_1(t)]\cup[h_2(t),1]\right\}\;.
\ee
In particular, ${\bf V}$ is of class ${\cal C}^{1,2}$.
\end{prop}

\begin{rem}{\rm The above representation  for ${\bf V}$ does not hold true in general.
For instance let $H(t,p)=\lambda(t)p(1-p)$ where $\lambda:[0,T]\to\R$ is Lipschitz continuous. We set $\Lambda(t)=\int_t^T\lambda(s)ds$. If
$$
\mbox{\rm $\lambda>0$ on $[0,b)$, $\lambda<0$ on $(b,T]$, $\Lambda(a)=0$}
$$
for some $0<a<b<T$, then one easily checks that
$$
{\bf V}(t,p)=\left\{\begin{array}{ll}
0 & {\rm if } \; t\in [0,a]\;, \\
\Lambda(t)p(1-p) & {\rm if } \; t\in [b,T]\; .
\end{array}\right.
$$
In particular
$$
{\bf V}(t,p)\neq \int_t^T{\rm Vex} H(s,p)ds =\Lambda (b)p(1-p) \qquad \forall (t,p)\in (a,b)\times (0,1) \; .
$$
Note also that in this example the dynamics is smooth, the value function ${\bf V}$ is smooth with respect to the variable $p$, but ${\bf V}$
just Lipschitz continuous with respect to the time variable.
}\end{rem}

\noindent \underline{Proof of Proposition \ref{ex1H}:} Let $w:[0,T]\times [0,1]\to \R$ be defined by
$$
w(t,p)=\int_t^T {\rm Vex} H(s,p)ds \qquad \forall (t,p)\in [0,T]\times \Delta(I)\;.
$$
We note that $w(T,p)=0$, $w(t,0)={\bf V}(t,0)$ and $w(t,1)={\bf V}(t,1)$.
One easily checks that $w$ is a solution of the HJ equation
$$
\left\{\begin{array}{l}
\min\left\{ w_t+H(t,p)\, , \, \frac{\partial^2 w}{\partial p^2}\right\}=0\\
w(T,p)=0
\end{array}\right.
$$
Indeed, if $p\in (h_1(t),h_2(t))$, then
$$
\frac{\partial^2 w}{\partial p^2}(t,p)=0 \; {\rm and }\; w_t(t,p)=H(t,h(t)) <H(t,p)\;.
$$
If $p\in (0,h_1(t)]\cup [h_2(t),1)$, then
$$
\frac{\partial^2 w}{\partial p^2}(t,p)\geq 0\; {\rm and }\; w_t(t,p)=H(t,p)
$$
(where the first equality holds in the viscosity sense since $w$ is convex with respect to $p$). Therefore $w={\bf V}$ and
${\cal H}$ is the set of points $(t,p)$ at which $H={\rm Vex}H$, i.e.,
given by (\ref{HpartCase}).
\qed

\begin{prop}\label{ex1p} Under the assumptions of Example \ref{ex1}, there is a unique optimal martingale measure $\bar \P$.
Under this martingale measure, the process $\p$ is purely discontinuous and satisfies:
$$
\p(t^-)=p_0 \qquad \forall t\in [t_0,t^*]\; \bar \P-\mbox{a.s., where }\; t^*=\inf\left\{t\geq t_0\; |\; p_0\in [h_1(t),h_2(t)]\right\}
$$
and
$$
\p(t)\in \{h_1(t),h_2(t)\} \qquad \forall t\in [t^*,T)\; \bar \P-\mbox{a.s.}
$$
In particular,
\be\label{proba}
\bar \P\left[ \p(t)=h_1(t)\; |\; \p(s)=h_1(s)\right] = \frac{h_2(t)-h_1(s)}{h_2(t)-h_1(t)} \qquad \forall t^*\leq s\leq t<T\;.
\ee
\end{prop}

\begin{rem}{\rm Set $T=\frac 14$, $h_1(t)=\frac 12-\sqrt t$, $h_2(t)=\frac 12+\sqrt t$, $t\in[0,T]$, $t_0=0$ and $p_0=\frac 12$. Since there is only one martingale measure that charges 
the graphs of $h_1$ and $h_2$, the process $\p$ under $\bar \P$ is,  up to a constant, the Az\'ema martingale with parameter 2 (see Emery \cite{emery}): 
under $\bar \P$, $(X(t):=\p(t)-\frac 12,\; t\in[0,T])$ satisfies the structure equation
\[ d[X]_t=dt-2X(t-)dX(t), \; t\in[0,T], \qquad X(0)=0.\]

}\end{rem}

\noindent \underline{Proof of Proposition \ref{ex1p}}: We do the proof in the case $t^*<T$ and $p_0\notin [h_1(t_0),h_2(t_0)]$, the proof of the other cases being similar. Under 
 these assumptions, $t^*>t_0$.
Let us fix $\bar \P$ some optimal martingale measure.
We need below the following result:

{\it Claim : } Let $\theta\geq t_0$ be a stopping time and let us assume that $\p(\theta^-)\notin[h_1(\theta),h_2(\theta)]$
on some set $A\in\FR_\theta$ with positive probability. Let $\theta'$ be  the stopping time
$$
\theta'=\inf\{t\geq \theta \; |\; \p(t)\in [h_1(t),h_2(t)] \ \}
$$
(by convention, $\theta'=T$ if there is no such a $t$). Then, on $A$,  $\bar \P-$a.s., $\theta'>\theta$ and $\p(t)=\p(\theta^-)$
for $t\in [\![\theta, \theta'[\![$.

{\it Proof of the claim : } Since $\p(\theta^-)\notin[h_1(\theta),h_2(\theta)]$ on $A$, ${\bf V}(\theta,\cdot)$ is strictly convex in a neighborhood of $\p(\theta^-)$.
Applying the equality obtained in Theorem \ref{bpH}:
\be\label{Eqtt-}
{\bf V}(t,\p(t))-{\bf V}(t,\p(t^-))-\langle \frac{\partial {\bf V}}{\partial p}(t,\p(t^-)), \p(t)-\p(t^-)\rangle =0
\ee
at $t=\theta$, we get that $\p(\theta)=\p(\theta^-)$ $\bar \P-$a.s. on $A$.
 Since $\p$ is right-continuous, $\theta'>\theta$ on $A$.
Using (\ref{Eqtt-}) again on $[\![\theta,\theta'[\![$ shows that $\p$ is continuous on $[\![\theta,\theta'[\![$, that $p(\theta'^-)\in\{ h_1(\theta'),h_2(\theta')\}$ (since $\p$ is not allowed to jump 
until it reaches the graphs of $h_1$ and $h_2$) and that, on $A$, we have 
$\theta'=\lim_{\epsilon\searrow 0}\theta'_\epsilon$, with, for $\theta'_\epsilon=\inf\{ t\geq \theta,\p(t)\in[ h_1(t)-\epsilon, h_2(t)+\epsilon]\}$, with $\theta'_\epsilon=T$ 
if there is no such $t$.\\
Let us now apply It\^{o}'s formula between $\theta$ and $\theta'$:
\be \label{AA}
{\bf E}_{\bar \P}\left[ {\bf V}(\theta',\p(\theta'))\right] =   {\bf E}_{\bar \P}\left[ {\bf V}(\theta,\p(\theta)) +\int_{\theta}^{\theta'}{\bf V}_t(s,\p(s))ds+
\frac12 \int_{\theta}^{\theta'} \frac{\partial^2 {\bf V}}{\partial p^2}(s,\p(s^-))d<\p^{c}>_s  \right]\;,
\ee
where $\p^c$ is the continuous part of $\p$ under $\bar \P$.
Since $\p(s)\in {\cal H}(s)$ for almost all $s$ $\bar P-$a.s.,
\be\label{BB}
{\bf E}_{\bar \P}\left[\int_{\theta}^{\theta'}{\bf V}_t(s,\p(s))ds\right]= - {\bf E}_{\bar \P}\left(\int_{\theta}^{\theta'}H(s,\p(s))ds\right)\;.
\ee
>From our assumption on ${\bf V}$, we also have
\be \label{DD}
 \frac{\partial^2 {\bf V}}{\partial p^2}(s,\p(s^-)) > 0 \qquad \forall s\in [\![\theta,\theta'[\![\; \bar P-{\rm a.s.,}
\ee
because
$$
 \frac{\partial^2 {\bf V}}{\partial p^2}(s,p)=\int_s^{t^*(p)} \frac{\partial^2 H}{\partial p^2}(\tau,p)d\tau  \qquad \mbox{\rm (where $t^*(p)=\inf\{s\geq 0\; |\; (s,p)\notin {\cal H}\}$),}
$$
 which is positive as soon as $t^*(p)>s$. Since, from Corollary \ref{PgrDynV}, we have
$$
{\bf E}_{\bar \P}\left[ {\bf V}(\theta',\p(\theta'))\right]={\bf E}_{\bar \P}\left[ {\bf V}(\theta,\p(\theta)) - \int_{\theta}^{\theta'}H(s,\p(s))ds\right]\;,
$$
combining the above equality with  (\ref{AA}), (\ref{BB}), (\ref{DD})  gives that $d<\p^{c}>_s=0$ $\bar P-$a.s. on $[\![\theta,\theta']\!]$.
This implies that,  for all $\epsilon>0$, the restriction on $A$ of the martingale $(\p_{(t\vee\theta)\wedge\theta'_\epsilon})_{t\in[0,T]}$ is simultaniously continuous and purely discontinuous
on $[0,T]$, thus it is constant. Therefore $\p$ restricted to $A$ is constant, equal to $\p(\theta^-)$ on $[\![\theta,\theta'[\![$, and the claim is proved. \\

Let us now prove that $\p(t)=p_0$ on $[t_0,t^*)$. For this we introduce the stopping time
$$
\theta=\inf\{ t\geq t_0 \; |\; \p(t)\in [h_1(t),h_2(t)]\ \}\;.
$$
Since $\p(t_0^-)=p_0\notin [h_1(t_0),h_2(t_0)]$, applying the claim to the stopping time $t_0$ we have that $\theta>t_0$
 and $\p(t)=p_0$ for $t\in [\![t_0,\theta[\![$ $\bar \P-$a.s. Since $\p(t)=p_0\notin [h_1(t),h_2(t)]$ for $t\in [\![t_0,\theta[\![$, we have $\theta\leq t^*$.
Since $\p(\theta)\in {\cal H}(\theta)$, we also have $\theta\geq t^*$.
Therefore $\theta=t^*$ and $\p(t)=p_0$ on $[t_0,t^*)$.\\

We now prove that $\p(t)\in \{h_1(t),h_2(t)\}$ for $t\in[t^*,T)$ $\bar \P-$a.s. Let us introduce, for any $\ep>0$,
the stopping time
\[ \theta_\epsilon=\inf\{ t\geq t^*|\p(t)\in [0,h_1(t)-\epsilon]\cup[h_2(t)+\epsilon,1]\},\]
(we set $\theta_\epsilon=T$ if there is no such a $t$).\\
Suppose now that there exists some $\epsilon>0$ such that $\bar \P[\theta_\epsilon<T]>0$. Without loss of generality we can suppose that the set $A:=\{ h_2(\theta_\epsilon)+\epsilon\leq\p(\theta_\epsilon)\leq 1\}\cap\{\theta_\epsilon<T\}$ satisfies $\bar P[A]>0$.
By definition of $\theta_\epsilon$, and since $h_2$ is increasing, we have
\begin{equation}
\label{infi}
 \p(t)<\p(\theta_\epsilon)\mbox{ on } [\![ t^*,\theta_\epsilon[\![\cap ([0,T]\times A).
\end{equation}
>From (\ref{Eqtt-}) again applied on $A$ at time $\theta_\epsilon$ we get, $\bar \P$-a.s., $\p(\theta_\epsilon^-)=\p(\theta_\epsilon)$ and therefore $\theta_{\tilde\epsilon}<\theta_\epsilon$, for all $0<\tilde\epsilon<\epsilon$.
But, still by (\ref{Eqtt-}) and the claim, on $A$, $\p$ is constant on the time interval $[\![\theta_{\tilde\epsilon},(\theta_{\tilde\epsilon})'[\![$.
Choosing now $\tilde\epsilon$ close enough to $\epsilon$ to get 
$0<\bar P[A\cap\{h_2(\theta_\epsilon)<p(\theta_{\tilde\epsilon})\}]\leq\bar P[A\cap\{ \theta_\epsilon<(\theta_{\tilde\epsilon})' \}]$, we obtain a contradiction to (\ref{infi}).\\

Equality (\ref{proba}) is then a straightforward application of the fact that $\p$ is a martingale which lives the union of the graphs of
$h_1$ and $h_2$. Finally, $\p$ is purely discontinuous since it has finite total variations.
\cq

\subsection{Examples in higher  dimensions}

Example \ref{ex1} can be extended to higher space dimensions. The interesting feature for $I\geq 3$ is that there are several optimal
martingale measures in general. They can be purely discontinuous, as in the two-dimensional case, but they can also be continuous.

\begin{eg} \label{ex2} We assume that there exists a smoothly evolving and increasing familly of smooth open convex subsets $(K(t))_{t\in [0,T]}$
of $\Delta(I)$, whose closure in contained in the interior of $\Delta(I)$, such that, for any $t\in [0,T]$, 
$$
\mbox{\rm ${\rm Vex}H(t,p)=H(t,p)$ $\; \Leftrightarrow\; $   $p\notin K(t)$ }\;,\qquad 
\mbox{\rm $H(t,\cdot)$ is affine on  $K(t)$}
$$
and
$$
\frac{\partial^2 H}{\partial p^2}(t,p) \; \mbox{\rm  is definite positive for $p\notin\overline{K(t)}$.}
$$
\end{eg}

Following the proof of Proposition \ref{ex1H} we get:

\begin{prop}\label{ex2H} Under the assumptions of Example \ref{ex2},
$$
{\bf V}(t,p)=\int_t^T {\rm Vex}H(s,p)ds \qquad \forall (t,p)\in [0,T]\times \Delta(I)
$$
and
$$
{\cal H}=\left\{(t,p)\in [0,T]\times \Delta(I)\; |\; p\notin K(t)\right\}\;.
$$
In particular, ${\bf V}$ is of class ${\cal C}^{1,2}$.
\end{prop}

Next we investigate the optimal martingale measures.

\begin{prop}\label{ex2p} Under the assumptions of Example \ref{ex2}, any optimal martingale measure  $\bar \P$ has the following structure:
$$
\p(t^-)=p_0 \qquad \forall t\in [t_0, t^*] \qquad {\rm and }\qquad
\p(t)\in \partial K(t) \qquad \forall t\geq t^*, \; \bar \P-a.s.,
$$
where $t^*=\sup\{t\geq t_0\; |\; p_0\notin K(t)\}$. Moreover, there is an optimal martingale measure under which $\p$  is purely
discontinuous. If, in addition, the family $(K(t))_{t\in [0, T]}$ has a positive minimal curvature and if $p_0\notin K(t_0)$, then there is also an optimal
martingale measure under which $\p$  is continuous.
\end{prop}

\begin{rem}{\rm An interesting case is when the evolving set $t\to\partial K(T-t)$ is moving according to its mean curvature (in $\Delta(I)$). Indeed, in this case, there is
an explicit formula for the martingale. If $p_0\in \partial K(t_0)$, then there is a optimal martingale measure under which the process
satisfies
$$
d\p(t)=\sqrt{2}\left(I-\nu(t,\p(t))\otimes\nu(t,\p(t))\right)dW_t
$$
where $I$ is the identity matrix of size $(I-1)$, $(W_t)$ is an $(I-1)-$dimensional Browian motion living in the hyperplane spanned by $\Delta(I)$ and
$\nu(t,p)$ denotes the unit outward normal to $K(t)$ at $p\in \partial K(t)$ (see \cite{BCQ} , \cite{SoTo}).
}\end{rem}

\dem The proof of the structure condition on the optimal martigale measures follows the same lines as the proof of Proposition \ref{ex1p} and
the existence of the purely discontinuous martingale measure can be  established as in Proposition \ref{CompDisc} because the set ${\cal H}$
satisfies (\ref{HLipschitz}).

Let us now check that there is a continuous optimal martingale measure. As usual we start our construction by building a
discrete time process $\p^n$. Let us fix $n$ large, $\tau=T/n$ the time-step and $t_k=k\tau$ for $k=0, \dots, n$.
For simplicity of notations we only build a process for an initial position such that $t_0=0$ and $p_0\in \partial K(0)$.
Let $Z$ be the set of $z=(z_i)\in \R^I$ such that $\sum_iz_i=0$. Note that  $Z=T_{\Delta(I)}(p)$ for any $p\in Int(\Delta(I))$. 
Let $f:Z\to Z$ be a Borel measurable process such that $|f(z)|=1$ and $\langle f(z),z\rangle=0$ for any $z\in Z$.
We denote by $\pi(t,p)$ a Borel measurable selection of the projection of $p$ onto the boundary of $K(t)$.

We now start the construction of $\p^n$. We set $\p^n_{-1}=p_0$. If $\p^n_k\in \partial K(t_k)$ is build, then we set
$$
q_1=\p^n_k+\lambda_1 f(\pi(t_k, \p^n_k)-\p^n_k) \; {\rm and }\; q_2=\p^n_k-\lambda_2 f(\pi(t_k, \p^n_k)-\p^n_k)
$$
where $\lambda_1>0$ and $\lambda_2>0$ are such that $q_1,q_2\in \partial K(t_{k+1})$.
Then we set $\p^n_{k+1}$ to be equal to $q_1$ with probability $\lambda_2/(\lambda_1+\lambda_2)$ and
$q_2$ with probability $\lambda_1/(\lambda_1+\lambda_2)$. The process $\p^n$ is then a martingale such that
$\p^n_k\in \partial K(t_k)$ for any $k=0, \dots n$.

We now show that, for any $\alpha>2$,  there is some constant $C_\alpha$ such that
\be\label{bebel}
{\bf E}\left[ |\p^n_{k_2}-\p^n_{k_1}|^\alpha\right] \leq C_\alpha |t_{k_2}-t_{k_1}|^{\alpha/2}\;.
\ee
Indeed from the Burkholder-Davis-Gundy inequality we have
$$
{\bf E}\left[ |\p^n_{k_2}-\p^n_{k_1}|^\alpha\right] \leq c_\alpha {\bf E}\left[ \left( <\p^n>_{k_2}-<\p^n>_{k_1}\right)^{\alpha/2}\right]
$$
$$
\leq c_\alpha {\bf E}\left[ \left(\sum_{k=k_1}^{k_2-1}|\p^n_{k+1}-\p^n_k|^2\right)^{\alpha/2}\right]
$$
Let us now estimate $|\p^n_{k+1}-\p^n_k|^2$. From the condition of positive minimal curvature, there is a constant $R>0$ such that
$$
K(t)\subset B_R(p-R\nu(t,p))\cap \Delta(I) \qquad   \forall p\in \partial K(t),\; \forall t\in [0,T]\;,
$$
where $\nu(t,p)$ is the outward unit normal to $K(t)$ at $p$.
Then, since $\p^n_{k+1}\in \partial K(t_{k+1})$, we have
$$
|\p^n_{k+1}- (\pi(t, \p^n_k)-R\nu(t_{k+1}, \pi(t,\p^n_k)))|\leq R\;,
$$
where
$$
\nu(t_{k+1}, \pi(\p^n_k))= \frac{\pi(\p^n_k)-\p^n_k}{|\pi(\p^n_k)-\p^n_k|}\;.
$$
Since by definition of $f$, $f(\pi(t,\p^n_k)-\p^n_k)\perp \pi(t,\p^n_k)-\p^n_k$, we have
$$
|\p^n_{k+1}-\p^n_k|^2+(R-|\pi(\p^n_k)-\p^n_k|)^2 \leq R^2\;,
$$
which implies that
$$
|\p^n_{k+1}-\p^n_k|^2\leq 2R L (t_{k+1}-t_k)
$$
where $L$ is a Lipschitz constant of the map $t\to K(t)$ in the Hausdorff distance. Therefore
$$
{\bf E}\left[ \left(\sum_{k=k_1}^{k_2-1}|\p^n_{k+1}-\p^n_k|^2\right)^{\alpha/2}\right]\leq  (2RL)^{\alpha/2}(t_{k_2}-t_{k_1})^{\alpha/2}\;.
$$
This proves (\ref{bebel}).

We now set $\p^n(t)=\p^n_k$ for $t\in[t_k, t_{k+1})$ and let $\P^n$ be the law of $\p^n$ on ${\bf D}(0)$. From Kolmogorov criterium we can extract a subsequence
of $(\P^n)$ which converges to some continuous martingale measure $\bar  \P$. Since $\p^n(t_k)\in \partial K(t_k)$ for any
$k=0, \dots, n$, we have $\p(t)\in \partial K(t)$ for any $t\in [0,T]$. Using It\^{o}'s formula and the fact that
$$
{\bf V}_t(s,p)=-H(s,p)
\;{\rm and }\;
\frac{\partial^2 {\bf V}}{\partial p^2}(s,p)=0 \qquad \forall p\in \partial  K(t), \; \forall t\in [0,T]\;,
$$
we have
$$
0={\bf E}_{\bar \P}\left[ {\bf V}(T,\p(T))\right]=
{\bf E}_{\bar \P}\left[ {\bf V}(0,p_0)+\int_0^T {\bf V}_t(s,\p(s))ds\right]
={\bf V}(0,p_0)- {\bf E}_{\bar \P}\left[ \int_0^T H(s,\p(s))ds\right]\;.
$$
Therefore
$$
{\bf V}(0,p_0)= {\bf E}_{\bar \P}\left[ \int_0^T H(s,\p(s))ds\right]\;,
$$
which shows that $\bar \P$ is optimal.
 \qed

\section{Conclusion}

In this paper we have investigated a continuous-time game with finite horizon and imperfect information on one side. We have proved that the
optimal behaviour of the informed player is directly related to the optimal revelation of his/her knowledge. This leads to an optimization problem
in which the unknown is a martingale measure. We have analysed this problem and found some necessary and some sufficient optimality conditions for
the optimal martingale measure.

Our analysis raises several intriguing questions: 
\begin{itemize}
\item We have seen that, under suitable regularity conditions, and, in particular, when $I=2$,
there are optimal martingale measures which are purely discontinuous. 
Does there always exist  an optimal martingale measure which is purely discontinuous~?

\item In our 2-dimensional examples, the optimal martingale measure is unique. Is this always the case when $I=2$ ?

\item In the case of a continuous time game in which both players have some private information,
existence and characterization  of the value are established in \cite{cr}. The equivalent of the martingale characterization (Theorem \ref{main})
in this framework is an open question.
\end{itemize}


\end{document}